\documentclass[journal]{IEEEtran}

\usepackage{ifpdf}

\usepackage{cite}

\ifCLASSINFOpdf
   \usepackage[pdftex]{graphicx}
\else
\fi

\usepackage[cmex10]{amsmath}

\usepackage{graphicx}
\usepackage{color}
\usepackage{subfig}
\usepackage{bm}
\usepackage{epstopdf}
\usepackage{epsfig}
\usepackage{float}
\usepackage{graphicx}
\usepackage{amsfonts,amssymb,amsmath}

\begin{document}

\title{Real-Time Monitoring of Area Angles with Synchrophasor Measurements}

\author{Wenyun~Ju,
Ian Dobson, 
Kenneth~Martin, Kai Sun, 
Neeraj Nayak, Iknoor Singh, \\Horacio Silva-Saravia, Anthony Faris, Lin Zhang,
Yajun Wang
\\
\vspace{-0.25cm}

\thanks{This work was supported by DOE under Award Number DE-OE0000849.

W. Ju, K. Martin, N. Nayak, I. Singh, and H.  Silva-Saravia are with Electric Power Group, LLC, Pasadena, CA USA (e-mail:wenyunju@gmail.com). 

I. Dobson is with Iowa State University, Ames, IA USA (e-mail:dobson@iastate.edu). 

K. Sun is with University of Tennessee, Knoxville, TN USA (e-mail:kaisun@utk.edu).  

A. Faris is with Bonneville Power Administration, Vancouver, WA USA (e-mail:ajfaris@bpa.gov). 

L. Zhang is with American Electric Power, New Albany, OH USA (e-mail:lin.zhang.ece@gmail.com). 

Y. Wang is with Dominion Energy, Richmond, VA USA (e-mail:yajun.wang@dominionenergy.com).
}
}
\maketitle
\begin{abstract}
This paper develops a comprehensive framework of Area Angle Monitoring (AAM) to monitor the stress of bulk power transfer across an area of a  power transmission system in real-time. Area angle is calculated from synchrophasor measurements to provide alert to system operators if the area angle exceeds pre-defined thresholds. This paper proposes general methods to identify these warning and emergency thresholds, and tests a mitigation strategy to relieve the area stress when the area angle exceeds the threshold. In order to handle the limited coverage of synchrophasor measurements, this paper proposes methods to estimate phase angles for boundary buses without synchrophasor measurements, which extends the application of AAM. AAM is verified for a power transmission area in the Western Electricity Coordinating Council system with both simulated data and synchrophasor measurements recorded from real events. A utility deployment to test the framework for monitoring area angle with live-stream and recorded synchrophasor data is described.
\end{abstract}

\begin{IEEEkeywords}
Real-time application, area angle, synchrophasor measurements, wide area monitoring, mitigation strategy 
\end{IEEEkeywords}

\IEEEpeerreviewmaketitle

\section{Introduction}
\IEEEPARstart{T}{he} stress of bulk power transfer through an area suddenly increases when there are multiple line outages inside the area.
Therefore it is important to monitor the stress in real-time and provide an alert to system operators if the outages cause overloads that make the system insecure. Then appropriate control actions can be taken to relieve the stress. 

Synchrophasor technology is developing rapidly in recent years. Synchrophasor technology uses monitoring devices, called phasor measurement units, which take high-speed measurements of phase angles, voltage and frequency that are time stamped with high-precision clocks. 
The measurements, typically taken 30 times a second, can quickly track system changes undetectable through traditional monitoring systems used in the industry \cite{FA1}. 
This makes new energy management applications possible, including model validation, oscillation monitoring, islanding detection and wide area monitoring.

References \cite{FA2,FA3} use angle difference between two buses to monitor the stress of power flow. 
Simulations of the system state before the 2003 USA/Canada blackout suggest the importance of increased angle difference for triggering blackouts \cite{FA2}. 
Indeed, many large cascading outages start with multiple outages initially occurring at a slow rate due to line overloads and other effects \cite{FA4,NERC02, 
AZFERCNERC12}.
Some cascading outages could be prevented if there are methods giving the situational awareness for fast emergency actions to relieve stress caused by multiple outages, and this is one motivation for real-time monitoring of area angles. 
The need for fast emergency actions can also arise when there are multiple simultaneous outages during extreme events such as storm, fire, icing, or an earthquake.

There are some established methods of detecting and resolving overloads due to power transfers.
References \cite{FA9,FA10,FA11} compute minimum security margins under operational uncertainty with respect to thermal overloads. Reference \cite{FA12} provides a tool for computation of transfer capability margins. 
These methods and applications are developed on top of SCADA and state estimation and can provide a comprehensive monitoring of the system status at the SCADA sampling rate.  Area Angle Monitoring (AAM) using synchrophasor measurements is approximate but faster, and it can monitor the stress of bulk power transfer through areas in real-time under multiple contingencies when the state estimator may not readily converge.

The concept of area angle is proposed in \cite{FA13,FA14} based on circuit theory. Reference \cite{FA15} shows that the area angle tracks the bulk power stress due to line outages inside the area and gives an approach to determine the emergency threshold of area angle. If the monitored area angle exceeds the emergency threshold, the area bulk power transfer should be reduced. 
Since reducing the area bulk power transfer can be translated into a specific action, the AAM is not only monitoring but also supplying actionable information to mitigate the stress.

Although some advantages of AAM have been demonstrated in \cite{FA15} using simulated data, there are still challenges for practically applying AAM in real-time using synchrophasor measurements. They are:

1) The deployment of AAM and testing using synchrophasor measurements in utilities need to be done.

2) AAM requires  synchrophasor measurements for all the boundary buses of the monitored area in order to calculate area angle. However, in reality, some boundary buses may not have synchrophasor measurements. 

3) The methods for quickly calculating area angle thresholds need improvements.

4) Mitigation strategies to reduce the area bulk power transfer need more detail and testing.

This paper addresses these challenges and moves AAM towards industry application. The main contributions are:

1) A comprehensive framework of AAM for real-time monitoring of bulk power transfer stress is proposed, which can guide utilities to apply AAM within their footprints. A real-time monitoring platform of AAM is developed and deployed at Bonneville Power Administration (BPA).

2) Either Linear State Estimation (LSE) or Phase Angle Compensation (PAC) is used to support the calculation of area angle when not all boundary buses have synchrophasor measurements, which extends the applicability of AAM. 

3) A new method to automatically calculate the warning threshold of area angle is proposed.

4) AAM is verified with both simulated data and synchrophasor measurements for an area in the Western Electricity Coordinating Council (WECC) system.

5) A practical mitigation strategy is proposed to release the stress of the bulk power transfer across the area when area angle exceeds the emergency threshold.

The rest of the paper is organized so that Sections II and III give an overview of AAM and a comprehensive framework of AAM, Section IV demonstrates AAM for an area in WECC, and Section V draws the conclusion.

\section{Overview of Area Angle Monitoring}
\begin{figure}[!htbp]
\centering
\includegraphics[width=45mm]{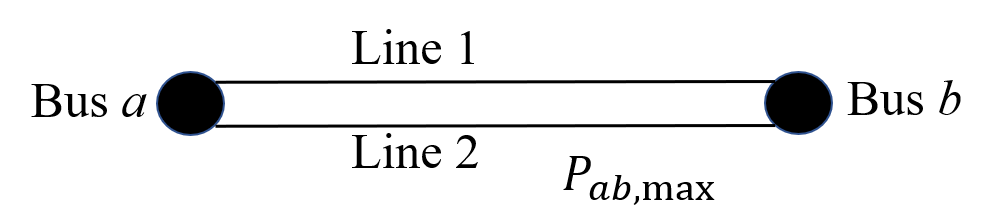}
\caption{Monitor stress across two parallel lines with angle difference.}
\label{Fig1}
\end{figure}

Two buses are connected by two identical lossless transmission lines in Fig.~\ref{Fig1}.
Assume $\theta_{ab}$ is the angle difference between buses $a$ and $b$, $P_{ab}$ is the power flow  from bus $a$ to bus $b$, and $P_{ab,max}$ is the maximum power flow from $a$ to $b$.

Consider two scenarios:

1. $P_{ab}$ increases.

2. Line 1 is tripped and $P_{ab}$ does not change. 

\noindent
For scenario 1, when $P_{ab}$ increases, $\theta_{ab}$ increases proportionally, indicating the increased stress of power flow. The maximum power flow  $P_{ab,max}$ 
does not change. 
In scenario 2, line 1 trips, $P_{ab,max}$ decreases and $\theta_{ab}$ increases. 
Therefore the angle difference $\theta_{ab}$ can indicate the increase of stress caused by either increased power flow or line outage. 
When a line outage occurs, $P_{ab,max}$ halves and $\theta_{ab}$ doubles. Thresholds for  $\theta_{ab}$ can be set up to distinguish outage severity.

\begin{figure}[!htbp]
\centering
\includegraphics[width=65mm]{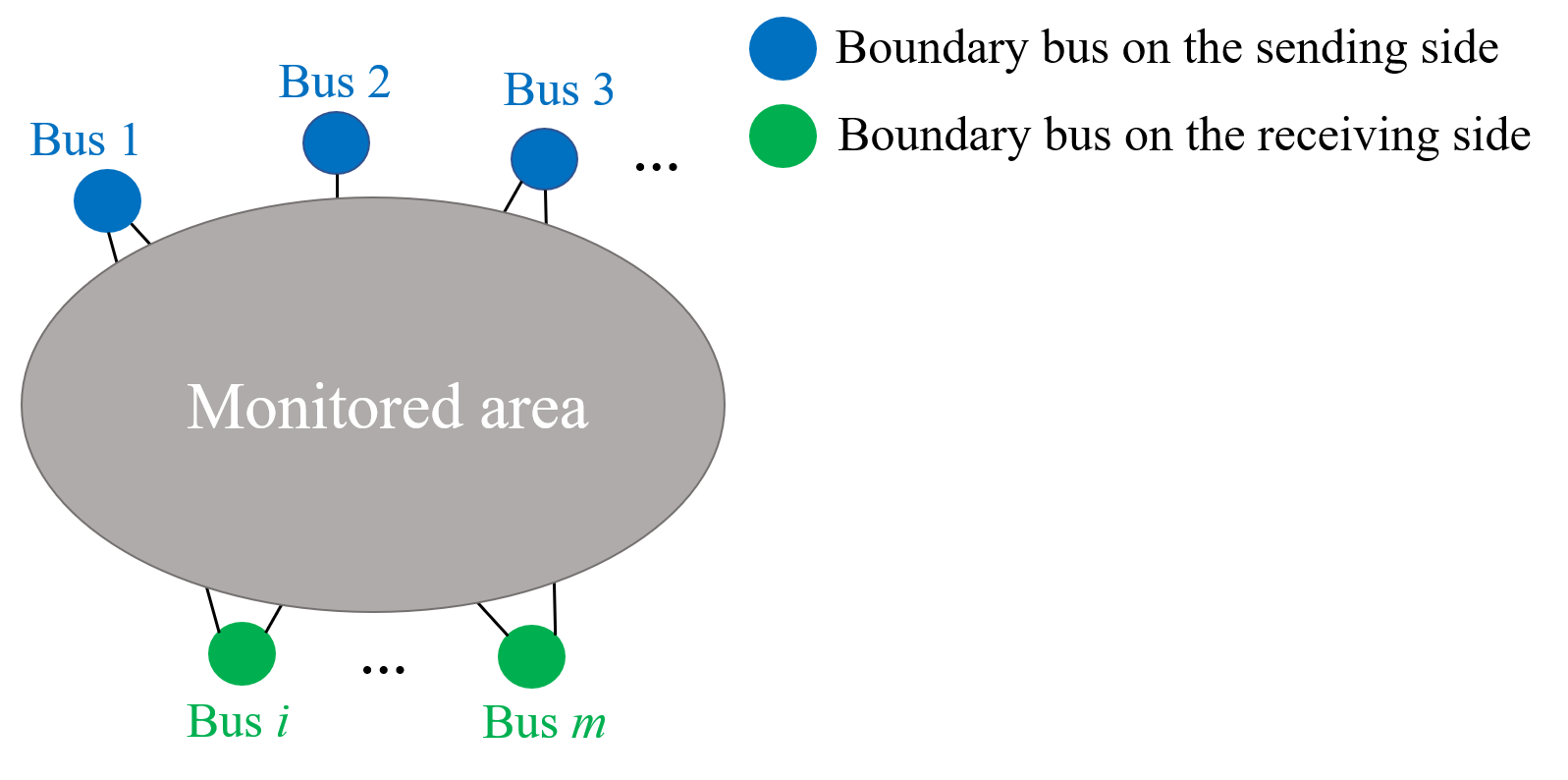}
\caption{Schematic diagram of a monitored area}
\label{Fig2}
\end{figure}

Area angle \cite{FA13,FA14,FA15} generalizes angle difference between buses to angle difference across an area as in Fig.~\ref{Fig2}. The stress of bulk power transfer through the area is indicated by a weighted combination of phase angles at the boundary buses of this area as the area angle:
\begin{eqnarray}
&\theta_{area} =
w_{1}\theta_{1}+w_{2}\theta_{2}+...+w_{m}\theta_{m}
\label{eq2}
\end{eqnarray}
where $w_{m}$ is the weight for bus $m$, $\theta_{m}$ is the phase angle of bus $m$, and $m$ is the number of boundary buses.

The weights on the boundary buses \cite {FA13} are calculated as
\begin{eqnarray}
w = (w_{1},w_{2},...w_{m}) =\frac{\sigma_{a}B_{eq}}{b_{mod}}
\label{eq3} 
\end{eqnarray}
where $\sigma_{a}$ is a vector with ones at the positions of the buses at the sending side of the area and zeros at the positions of the buses at the receiving side, $B_{eq}$ is the equivalent susceptance matrix of the area at the boundary buses, and $b_{mod} = \sigma_{a}B_{eq}\sigma_{a}^{T}$ is the bulk susceptance of the area. 

It is approximately the case that the monitored area angle gets larger as the maximum power that could enter the area decreases. This property can be used to set up alarm/warning and emergency thresholds to monitor the area stress  \cite{FA15}.

The reason for using the area angle to monitor stress for only one particular pattern of power flow through only one specific area is that if the area angle indicates too much stress, the mitigating action is clearly to reduce that particular power flow through the area \cite{FA15}. 
That is, monitoring the area angle for a specific area gives actionable information. 
Other area angles can be set up and monitored for other areas and patterns of stress as needed. In contrast to the area angle, more generic combinations or selections of angle differences can indicate a general stress but are not well associated with specific actions.

\section{Framework of Area Angle Monitoring}
\begin{figure}[!htbp]
\centering
\includegraphics[width=90mm]{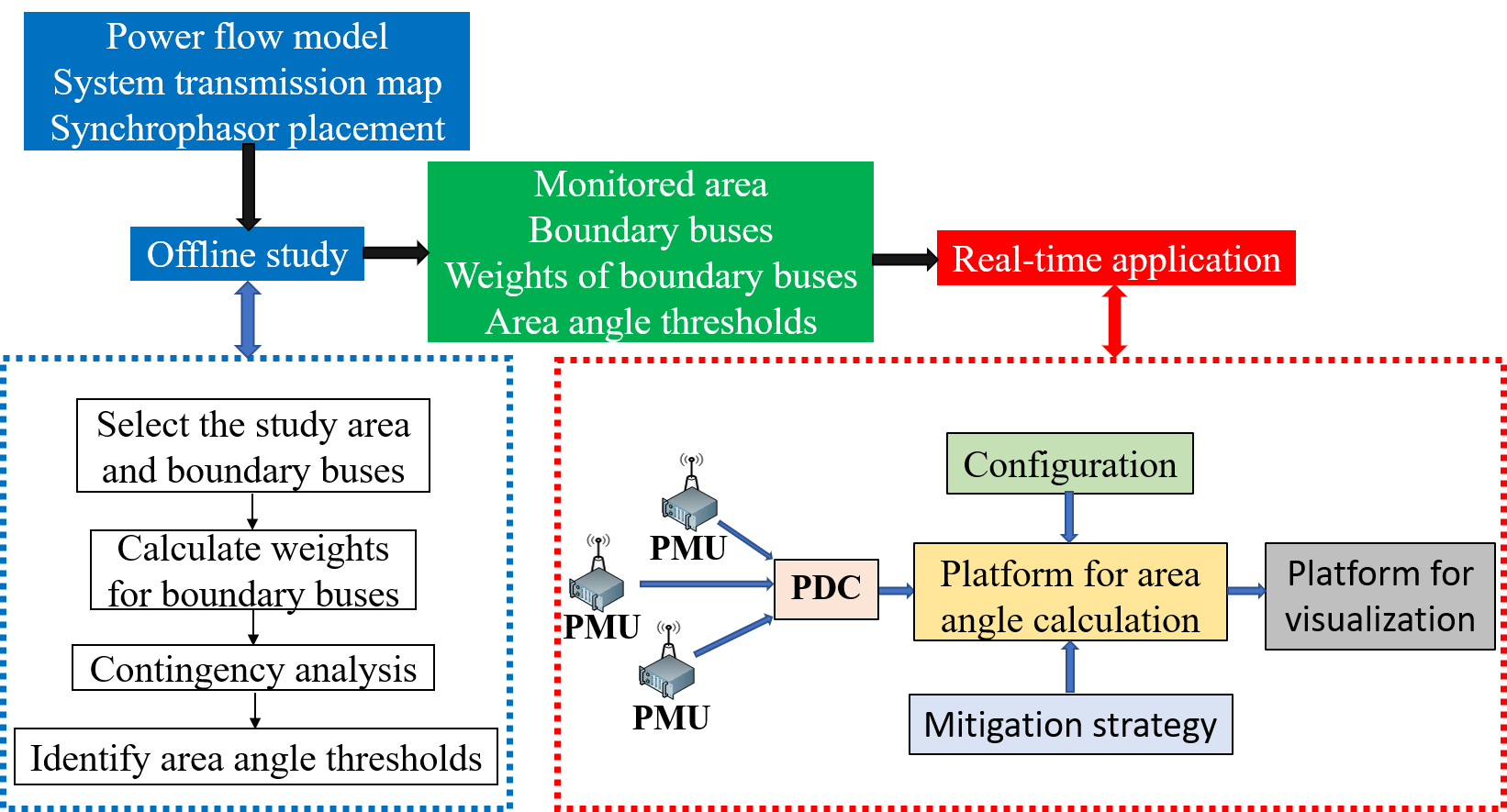}
\caption{Framework of AAM}
\label{Fig3}
\end{figure}

We propose the comprehensive framework of AAM shown by Fig.~\ref{Fig3}. It gives the core components (offline calculations and real-time application) needed for utilities to apply AAM in real-time  using synchrophasor measurements.

The offline study provides the monitored area, the boundary buses, weights of boundary buses, and area angle thresholds. They are needed for  calculation of area angle and detection of system status in real-time application. 
The inputs of the offline study include the power flow model, the system transmission map, and the synchrophasor placement. 

The offline study contains 4 steps:

%

%

%

Step 1: With the power flow model, select a monitored area, its sending and receiving edges, and its  boundary buses that mostly have synchrophasor measurements.

Step 2: Calculate the weights for the boundary buses with Kron reduction.

Step 3: Contingency analysis. Calculate the maximum powers that could enter the monitored area and the corresponding area angles under N-1 contingencies. For more details see \cite[Section III-D]{FA15}.

Step 4: Identify area angle thresholds based on the results from Step 3. This will be discussed in Section III-A. 


For real-time application, synchrophasor measurements are collected in a Phasor Data Concentrator (PDC) from PMUs in substations. 
These data are sent to a platform for calculating area angle. 
The calculation needs a pre-defined configuration, including the boundary buses, weights of boundary buses, and area angle thresholds. 
Then the results are visualized and shown to system operators. 
Warning or emergency status will be detected if the area angle exceeds the corresponding threshold. 
A mitigation strategy will be recommended to reduce the area stress if emergency status is detected. A real-time monitoring platform of AAM  is implemented and deployed at BPA, as described in Section IV-H.

Some low-quality synchrophasor measurements could affect the accuracy of the area angle.
However, low-quality synchrophasor measurements can be detected by the techniques in \cite{FA60} that can be easily integrated into the proposed framework. 
\vspace{-0.4cm}
\subsection{Area Angle Thresholds}
Reference \cite {FA15} defines an emergency threshold of area angle, but  the warning threshold is subjective. This paper proposes a better method to identify the warning threshold.

Consider the set of maximum powers that could enter the monitored area without violating line flow limits under \hbox{N-1} with $n$ line contingencies as $\{P_{mod}^{1}, ..., P_{mod}^{i}, ..., P_{mod}^{n}\}$ sorted into a descending order. 
(These maximum powers are calculated according to   \cite[Section III-D]{FA15}.) 
$P_{mod}^{1}$ corresponds to the least severe contingency, $P_{mod}^{i}$ corresponds to contingency $i$ and $P_{mod}^{n}$ corresponds to the most severe contingency. 

For contingency $i$, the system is placed in the limiting power flow condition of $P_{mod}^{n}$, line $i$ is outaged, and the corresponding area angle $\theta_{mod}^{i}$ is then calculated as:
\begin{eqnarray}
&\theta_{mod}^{i} =
w_{1}\theta_{1}^{i}+w_{2}\theta_{2}^{i}+...+w_{m}\theta_{m}^{i}
\label{eq4}
\end{eqnarray}
where $\theta_{m}^{i}$ is the phase angle of boundary bus $m$ under contingency $i$.
Doing this calculation for each contingency gives the set of area angles $\{\theta_{mod}^{1}, ..., \theta_{mod}^{i}, ..., \theta_{mod}^{n}\}$.

The standard deviation for three consecutive points in $\{P_{mod}^{1}, ..., P_{mod}^{i}, ..., P_{mod}^{n}\}$ is calculated as
\begin{align}
\sigma_{k}  &= \sigma([P_{mod}^{k-2}, P_{mod}^{k-1}, P_{mod}^{k}]),\quad  3\leq k\leq n
\label{eq40} 
\end{align}
starting with 
$k$=3 and increasing $k$ until $\sigma_{k}\geq \tau$, where $\tau$ is a constant. Then 
the warning threshold is
%
\begin{eqnarray}
&& \theta^{thr,w}_{mod} = \theta_{mod}^{k}
\label{eq30}
\end{eqnarray}
$\theta^{thr,w}_{mod}$ is the first point at which the maximum powers decrease significantly, indicating a relatively heavy stress inside the monitored area. Any other contingency that causes the area angle to be larger than $\theta^{thr,w}_{mod}$ will give the warning status. Also practically useful is that no action is needed if the area angle is less than $\theta^{thr,w}_{mod}$.

The area angle corresponding to $P_{mod}^{n}$  is identified as the emergency threshold \cite{FA15}:
\begin{eqnarray}
&& \theta^{thr,e}_{mod} = \theta_{mod}^{n}
\label{eq6}
\end{eqnarray}
$\theta^{thr,e}_{mod}$ corresponds to the largest area stress  satisfying the \hbox{N-1} security criterion. Any multiple contingencies that cause the area angle to be larger than $\theta^{thr,e}_{mod}$ will give emergency status since they correspond to violating the N-1 criterion.
\vspace{-0.3cm}
\subsection{Area Angle Thresholds with Angle Compensation}

The thresholds of area angle obtained from Section III-A. are based on the power flow model. There may be a mismatch of area angle calculated from the power flow model and real-time monitoring with synchrophasor measurements for normal status. The mismatch is defined as
\begin{eqnarray}
&& \Delta \theta_{com} = \theta_{ope} -\theta_{mod}
\label{eq7}
\end{eqnarray}
where $\theta_{mod}$ is the area angle obtained from the power flow model for normal status and $\theta_{ope}$ is the area angle obtained from real-time operation for normal status. 

By compensating the angle thresholds obtained from the power flow model with $\Delta \theta_{com}$, the area angle thresholds for real-time monitoring are calculated:
\begin{align}
 \theta^{thr,w}_{ope} &= \theta^{thr,w}_{mod} +\Delta \theta_{com}
\label{eq8}\\
 \theta^{thr,e}_{ope} &= \theta^{thr,e}_{mod} +\Delta \theta_{com} 
\label{eq9}
\end{align}

\subsection{Detection of Warning and Emergency Status}
For real-time monitoring, when the area angle $\theta_{area}$ exceeds the emergency threshold $\theta^{thr,e}_{ope}$, an emergency status is detected. 
When $\theta^{thr,w}_{ope}\leq\theta_{area}<\theta^{thr,e}_{ope}$, a warning status is detected. A time delay $t_{area}$ is applied to prevent false detection of warning or emergency status.

\subsection{Mitigation Strategy for Reducing the Stress}
If the area angle exceeds the emergency threshold, indicating that the stress of bulk power across the area violates the N-1 criterion, the stress needs to be mitigated quickly.
The area angle has an advantage of a physical interpretation as the angle across the area satisfying Ohm's law \cite{FA13}.
Ohm's law ensures that a mitigation strategy reducing the power flow through the area will reduce the area angle proportionally.
Another advantage of real-time area angle monitoring is that if operators perform the mitigation, they can quickly see the response of area angle to verify the mitigation.

One simple mitigation strategy is  generator ramp up or load shedding on the receiving side of the power system. This is equivalent to and can be tested by shedding load on the boundary buses at the receiving side of the area. 
Assuming the total amount of load to shed is $L_{total}$, we shed  load on each receiving bus is proportional to the magnitude of its weight:
\begin{eqnarray}
L_{j} = |w_{j}|L_{total}, \quad  1\leq j\leq r
\label{eq24}
\end{eqnarray}
where $r$ is the number of boundary buses at the receiving side. (Note that the boundary buses on the receiving side have negative weights that sum to $-1$.)

\subsection{Estimate Angle on Boundary Buses without Synchrophasors}

For the calculation of area angle, the ideal situation is that all boundary buses are installed with synchrophasor measurements. However, in reality, it is common that not all boundary buses are installed with synchrophasor measurements. 

\begin{figure}[!htbp]
\vspace{-0.3cm}
\centering
\includegraphics[width=60mm]{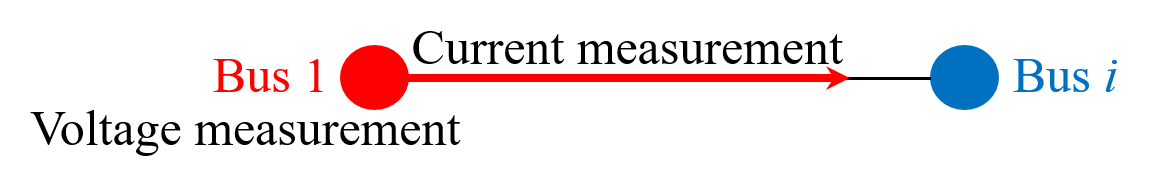}
\caption{Linear State Estimation of bus $i$ phase angle using voltage and current PMU measurements at neighboring bus 1.}
\label{Fig41}
\vspace{-0.3cm}
\end{figure}

Linear state estimation \cite {FA17} can extend the observability of system using synchrophasors
in the scenario shown by Fig.~\ref{Fig41}. In Fig.~\ref{Fig41}, bus $i$ is a boundary bus without synchrophasor measurements connected by a transmission line to bus 1 with synchrophasor measurements. 
The phase angle of bus $i$ can be estimated if bus 1 has voltage measurement and a current measurement on the transmission line and the impedance of the transmission line is known. 

\begin{figure}[!htbp]
\centering
\includegraphics[width=65mm]{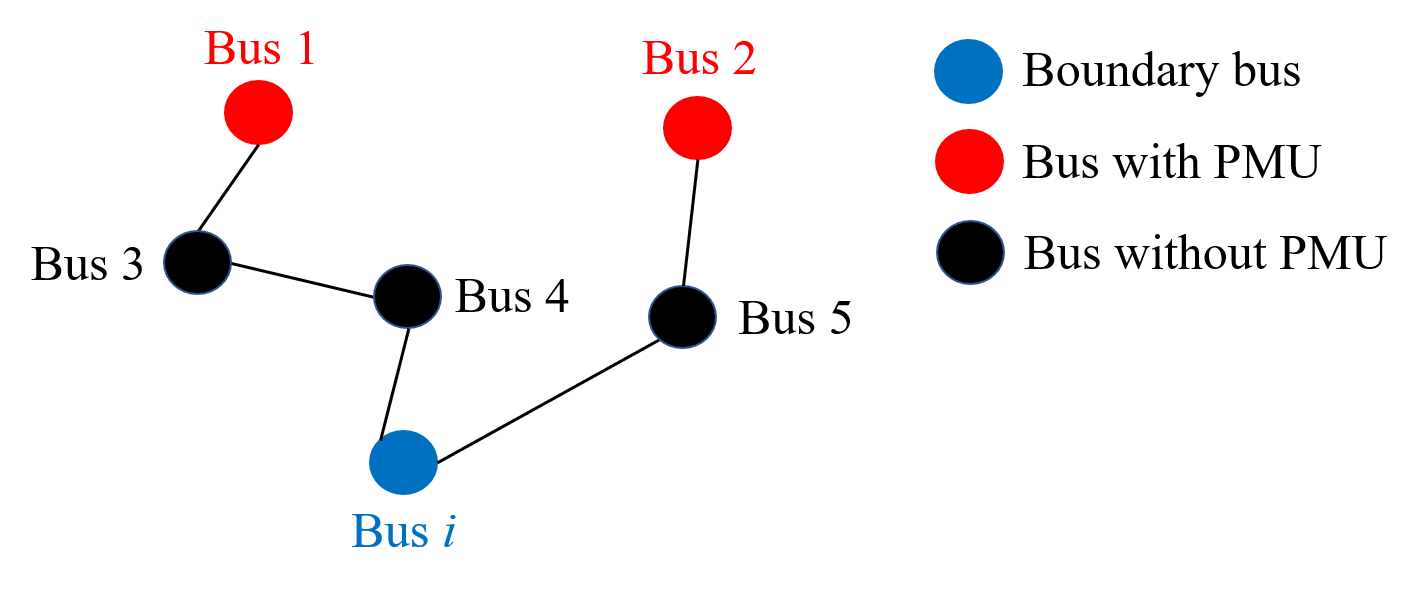}
\caption{Estimating phase angle of bus $i$ using PAC}
\label{Fig4}
\end{figure}

Consider the more general scenario shown by Fig.~\ref{Fig4}: buses 1 and 2 are far away from the boundary bus $i$ and they have synchrophasor measurements, and buses 3, 4, 5 without synchrophasor measurements are the neighbors of bus~$i$. The phase angle of bus~$i$ cannot be estimated using LSE.
However, a practical method to roughly estimate the phase angle of bus~$i$ using Phase Angle Compensation (PAC) is
\begin{eqnarray}
&&\theta_{i} = \theta_{j,sm} + \theta_{i,PAC}
\label{eq10}
\end{eqnarray}
where $\theta_{j,sm}$ is the phase angle of bus $j$ with synchrophasor measurements ($j$ is 1 or 2 in Fig.~\ref{Fig4}). 
The PAC $\theta_{i,PAC}$ of boundary bus $i$ is the angle difference between bus $i$ and bus~$j$ calculated from DC power flow in the offline study. 
Bus~$j$ is the  bus with synchrophasor measurements closest to bus $i$ in terms of electric distance.

The set of PACs is obtained from offline study and used in real-time, and we prefer to use a constant set of PACs rather than updating it frequently since it involves several steps from the offline study. 
The case studies in subsection~\ref{PACinfluence} test the accuracy of PAC for estimating area angle.

\subsection{Updating Thresholds when Topology Changes Significantly}
The area angle thresholds are computed offline and used in real-time. They need to be updated when the system topology changes significantly, such as in scheduled maintenance.

\begin{figure}[!htbp]
 \centering
\subfloat[Original method]{\includegraphics[width=36mm]{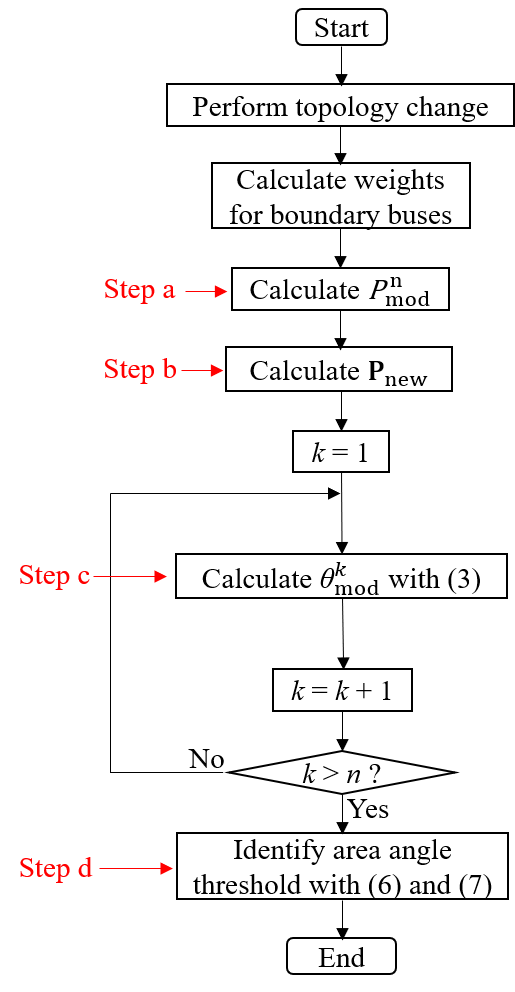}} \hspace{6mm}
\subfloat[Proposed method]{\includegraphics[width=29.5mm]{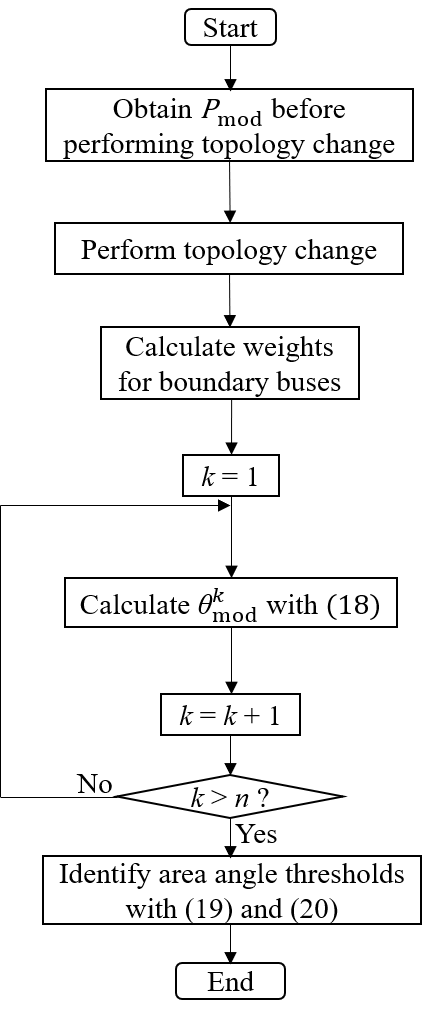}} 
\caption{Methods for updating area angle thresholds under topology change}
\label{Fig5}
\vspace{-0.3cm}
\end{figure}

\textit{1) Original Method}

Reference \cite{FA15} does not explicitly present a method to update area angle thresholds, but a similar approach can be summarized by Fig.~\ref{Fig5}(a). There are 4 important steps:

Step~$a$: Calculate the set of maximum powers $\{P_{mod}^{1}, P_{mod}^{2}, ..., P_{mod}^{n}\}$ for all N-1 contingencies sorted into a descending order. Select the maximum power $P_{mod}^{n}$ with the worst case contingency.

Step~$b$: Obtain a new bus injection vector $P_{new}$ by placing the system in the condition of limit of $P_{mod}^{n}$.

Step~$c$: For contingency $k$, calculate the area angle $\theta^{k}_{mod}$.

Step~$d$: Identify area angle thresholds with (\ref{eq30}) and (\ref{eq6}). Note that (\ref{eq30}) is not used by \cite{FA15}.

The area angle for contingency $k$ using the weights before contingency $k$ is calculated as (\ref{eq11}):
\begin{align}
 \theta^{k}_{mod}=w\theta^{k}_{m}=\frac{\sigma_{a}B_{eq}}{b_{mod}}\theta^{k}_{m}
\label{eq11}\\
 \theta^{k}_{m}=e_{ab}\cdot[(B^{k})^{-1}P_{new}]
\label{eq12}
\end{align}
where $\theta^{k}_{m}$ is the vector of phase angles for the boundary buses, $B^{k}$ is the susceptance matrix under contingency $k$, $e_{ab}$ is the column vector of length $n$ with ones at the positions for the boundary buses.

\textit{2) Proposed Method}

The original method needs steps a,b,c,d to calculate the area angle for each contingency. Among these, step a is especially complex and it requires the calculation of $\{P_{mod}^{1}, P_{mod}^{2}, ..., P_{mod}^{n}\}$. In order to avoid this calculation, we propose the simplified approximate method of Fig.~\ref{Fig5}(b).

If the topology change caused by contingency $k$ is considered in the calculation of weights, the calculated area angle using the updated weights is indicated by $\theta^{[k]}_{mod}$, where
\begin{eqnarray}
\theta^{[k]}_{mod} = \frac{\sigma_{a}B_{eq}^{k}}{b_{mod}^{k}}\theta^{k}_{m}
\label{eq13}
\end{eqnarray}
$B_{eq}^{k}$ is the equivalent susceptance matrix of the boundary buses considering contingency $k$, and $b_{mod}^{k}$ is the bulk susceptance of the area considering contingency $k$.

Apply Ohm's law to the area angle at the maximum power transfer of the base case to get 
\begin{eqnarray}
&& P_{mod}=b_{mod}\theta_{mod}
\label{eq14}
\end{eqnarray}
where $P_{mod}$ is the base case maximum power through the area without violating the line flow limits, and $\theta_{mod}$ is the corresponding area angle with $P_{mod}$ through the area.

Also, considering the maximum power transfer when contingency $k$ occurs, we have
\begin{eqnarray}
& P_{mod}^{k}=b_{mod}^{k}\theta_{mod}^{[k]}
\label{eq15}
\end{eqnarray}
where $P_{mod}^{k}$ is the maximum power through the area considering contingency $k$, and $\theta_{mod}^{[k]}$ is the corresponding area angle with $P_{mod}^{k}$ through the area.

The results from section 4.4 in \cite{FA18} show that $\theta^{[k]}_{mod} \approx \theta^{k}_{mod}$. Then we can have:
\begin{eqnarray}
&& P_{mod}^{k}=b_{mod}^{k}\theta_{mod}^{k}
\label{eq20}
\end{eqnarray}
Sometimes there are no parallel paths outside the area for power transfer from the sending side to the receiving side of the area. The area is a cutset area \cite{FA13,FA14}.
Then when a non-islanding contingency $k$ occurs, the generation stays the same and $P_{mod}^{k}=P_{mod}$.
Moreover, there is a good approximation when the parallel paths outside the area have high impedance. In this case, when a non-islanding contingency $k$ occurs, $P_{mod}^{k}$ $\approx$ $P_{mod}$  \cite {FA18,FA19}. So we have:
\begin{eqnarray}
&& \theta_{mod}^{k} = 
\frac{P_{mod}^{k}}{b_{mod}^{k}} \approx \frac{P_{mod}}{b_{mod}^{k}}
\label{eq16}
\end{eqnarray}
This approximation, but applied to the monitored area angle instead of the threshold angle, is also used and discussed in \cite[eqn. (13)] {FA19}.

Instead of using (\ref{eq30}) and (\ref{eq6}), we propose the approximate emergency threshold $\hat{\theta}_{mod}^{thr,e}$ and warning threshold $\hat{\theta}_{mod}^{thr,w}$ as
%
\begin{align}
 \hat{\theta}_{mod}^{thr,e} &= \max\{\theta_{mod}^{1},\theta_{mod}^{2}, ... ,\theta_{mod}^{n}\} 
\label{eq17}\\
 \hat{\theta}_{mod}^{thr,w}& = \frac{1}{2}\big[\hat{\theta}_{mod}^{thr,e}
 +\min\{\theta_{mod}^{1},
 ...,\theta_{mod}^{n}\}\big]
\label{eq18}
\end{align}
The warning threshold (\ref{eq18}) is
the average value of the minimum and maximum values of $\{\theta_{mod}^{1},\theta_{mod}^{2}, ... ,\theta_{mod}^{n}\}$.
Since the emergency threshold is used to detect more severe events, its approximation (\ref{eq17}) is acceptable only if it is sufficiently accurate. The case study in subsection \ref{testtopchange} compares the approximate and original area angle thresholds.

\section{Case Studies}
\subsection{Model and Parameter Preparation}
The case studies use the power flow model for the 2020 heavy summer case of the WECC system provided by BPA.

PMUs are mainly deployed for high-voltage power transmission. For a reduced power flow model with high-voltage level, it is relatively easy to select a monitored area with a large coverage of synchrophasor measurements on the boundary buses. However, the model we have is a detailed model, and under this circumstance, static network reduction is needed.

A reduced/equivalent model ($\geq$ 230kV) is obtained using the modified Ward reduction \cite{FA16}. Note that Kron reduction is a standard tool to obtain a ``network-reduced" or ``Ward-equivalent" model for power flow studies \cite{FA27}. The main difference between the Kron/Ward reduction and the modified Ward reduction is that all generators in the original model are retained integrally in the reduced model with the modified Ward reduction.
Reference \cite{FA16} verifies the accuracy of the modified Ward reduction by comparing the power flows of the original and reduced models.
The original and reduced models are compared in Table~\ref{tab1}.

\begin{table}[!htbp]
\caption{Comparison of Models}
\centering
\begin{tabular}{ccc}
& Original Model& Reduced Model  \\ \hline 
Number of Buses& 20507& 3101 \\ 
Number of Generators&4019& 4014 \\  
Number of Lines & 26395 & 8000  \\ \hline
\end{tabular}
\label{tab1}
\vspace{-0.3cm}
\end{table}

The parameters chosen are 
$\tau=0.5$ p.u. and $t_{area}= 5$\,s. The offline study of AAM in Fig.~\ref{Fig3} is implemented with MATLAB R2019a. GE Positive Sequence Load Flow Software V21.5 is used to generate the simulated data used in subsections IV-C and IV-E.  Light, medium, and heavy loadings are considered. The difference of active load between heavy and medium loadings is 1000 MW, and the difference of active load between medium and light loadings is around 23\,000 MW.

\vspace{-0.1cm}
\subsection{Monitored Area and Area Angle Thresholds}

A monitored area is selected inside the reduced model of WECC. The monitored area roughly covers Oregon state  and contains 176 transmission lines and 106 buses. The bulk power transfer of interest is from north (sending side) to south (receiving side). There are 14 boundary buses; 7 of them are on the sending side and 7 are on  the receiving side.
The area angle weights of buses 1--7 on the sending side are [\text{0.1271}, \text{0.5303}, \text{0.2616}, \text{0.0396}, \text{0.0385}, \text{0.0005}, \text{0.0023}], and the weights of buses 8--14 on the receiving side are [\text{--0.1269}, \text{--0.0958}, \text{--0.0017}, \text{--0.1615}, \text{--0.2979}, \text{--0.2766}, \text{--0.0395}].

\textit{1) Max Powers and Area Angles under N-1 Contingencies}

\begin{figure}[!htbp]
\centering
\vspace{-0.5cm}
\begin{tabular}{c c}
\subfloat[N-1 contingencies; \newline equivalenced lines excluded ]{\includegraphics[width=44mm]{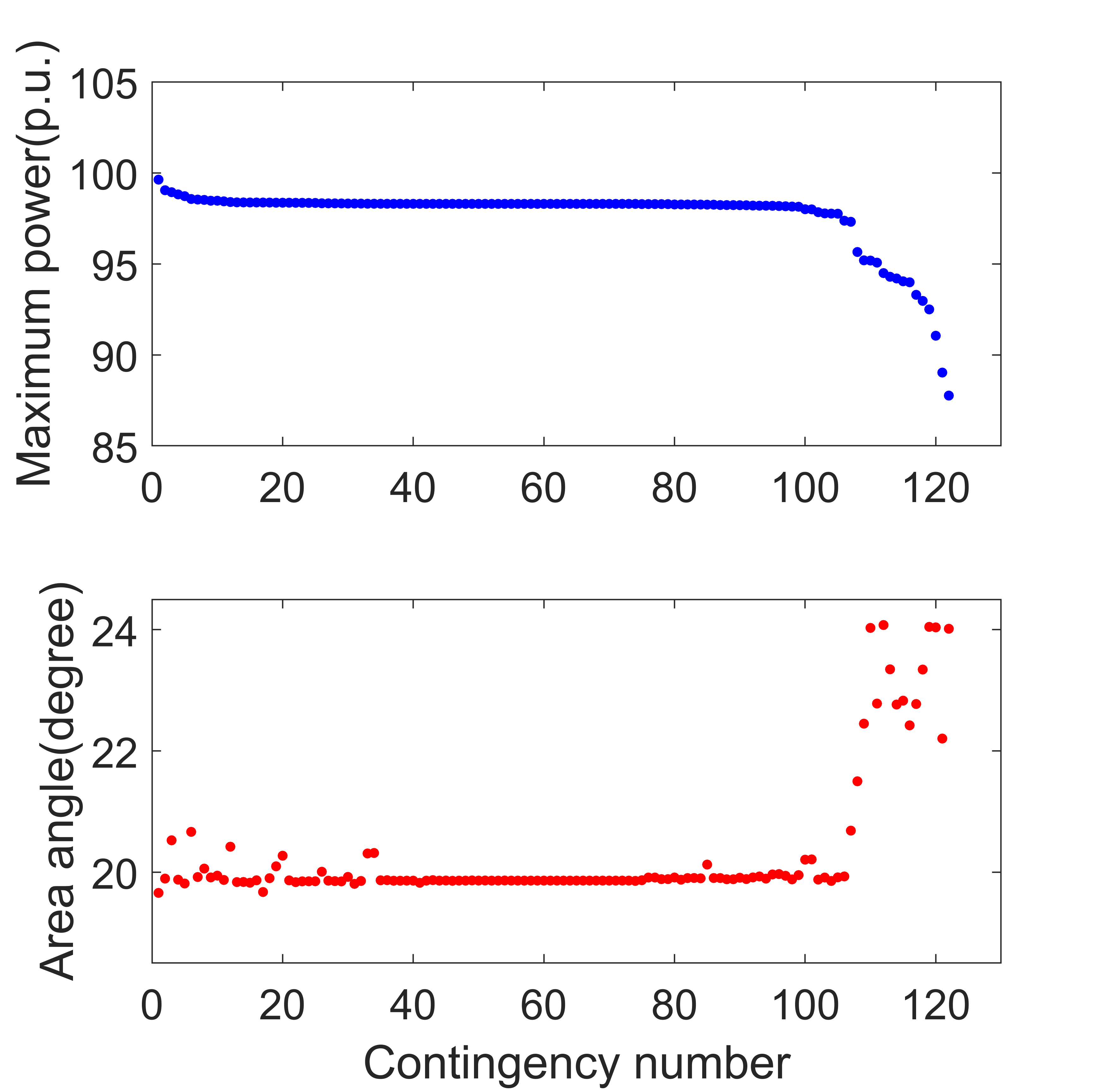}} \
\subfloat[N-1 contingencies; \newline all lines in reduced model]{\includegraphics[width=44mm]{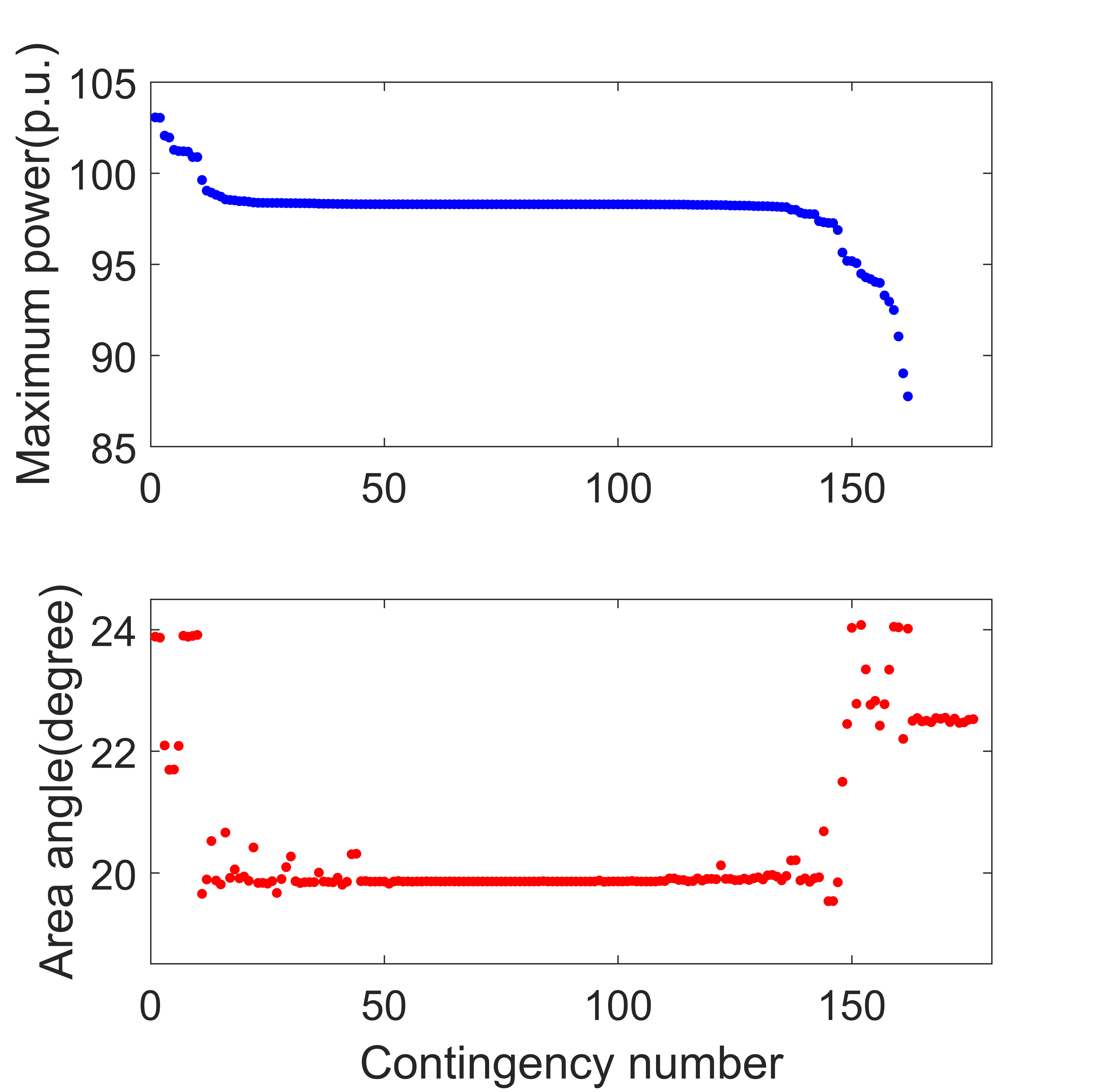}} 
\centering
\end{tabular}
\caption{Maximum powers and area angles}
\label{Fig24}
\vspace{-0.3cm}
\end{figure}

To set the emergency threshold, we need to examine the worst case maximum powers that could enter the monitored area under N-1 line contingencies.
The maximum powers that could enter the monitored area and area angles corresponding to the non-islanding N-1 contingencies are shown by Fig.~\ref{Fig24}(a). Basically the area angle increases as the maximum power decreases \cite{FA15}. This verifies that area angle can distinguish the stress of bulk power transfer caused by different contingencies. Fig.~\ref{Fig24}(a) is used to set the emergency threshold.

\looseness=-1
Note that Fig.~\ref{Fig24}(a) excludes lines that are equivalenced in the reduced model; that is, we only apply the N-1 criterion to lines within the monitored area in the reduced model that also appear in the original model. 
The reason is that removing an equivalenced line in the reduced model has an effect that is unrelated to  the effect of removing a line in the original model, so that applying the N-1 criterion with the equivalenced lines does not correctly reflect the N-1 criterion applied to the real system.
We can see the effect of applying the N-1 to all the lines within the monitored area of the reduced system, including the 54 equivalent lines (176 lines inside the area), in Fig.~\ref{Fig24}(b), as additional more extreme outliers.
(Of course, one way to prevent problems with equivalenced lines is to avoid system reduction, but that entails a larger system model.)

We check that the 6581 lines of the detailed model eliminated in the system reduction do not significantly affect the area angle by 
removing each of those lines in the detailed model, obtaining a new reduced model with each of those lines in the detailed model removed, and recalculating the area angle  with the system placed in the condition of limit of $P_{mod}^{n}$.
The recalculated area angles are within 0.4 degree of the baseline of area angles in Fig.~\ref{Fig24}(a), indicating that the \hbox{N-1} contingencies of the lines eliminated in the reduction have little effect on the stress inside the monitored area.


\textit{2) Area Angle Thresholds}

The warning threshold using (\ref{eq40})-(\ref{eq30}) and emergency threshold using (\ref{eq6}) are calculated as  21.49 degree and  24.07 degree. 

The value of $\Delta\theta_{com}$ is calculated as --0.7 degree, and is used to adjust the warning and emergency thresholds for real-time monitoring with (\ref{eq8}) and (\ref{eq9}) as shown in Table~\ref{tab2}. In the following sections, $\theta^{thr,w}_{ope}$ and $\theta^{thr,e}_{ope}$ are used for the warning and emergency thresholds.

\begin{table}[!htbp]
\centering  
\vspace{-0.2cm}
\caption{Area Angle Thresholds}
\begin{tabular}{cc}
Threshold & Value of Threshold (degree) \\ \hline 
$(\theta^{thr,w}_{mod}, \theta^{thr,e}_{mod})$& (21.49, 24.07) \\ 
$(\theta^{thr,w}_{ope}, \theta^{thr,e}_{ope})$& (20.79, 23.37)  \\   \hline
\end{tabular}
\label{tab2}
\vspace{-0.3cm}
\end{table}

\vspace{-0.1cm}
\subsection{Verification of AAM With Simulated Data}
AAM is verified with simulated data. Phase angles obtained from dynamic simulation can be used as fictitious synchrophasor measurements. The method in Section III-E is used.
This subsection considers  contingencies both outside and inside the monitored area that  stress  the area. 

\textit{1) Verification of AAM with Generator Trip}

Contingency 1: Trip of one generator (around 1400\,MW output) in the southern part of the receiving side of the monitored area at 60\,s.

Contingency 2: Trip of two generators (around 2800\,MW output) in the southern part of the receiving side of the monitored area at 60\,s.

These contingencies are outside the area and they reduce the generation on the receiving side of the monitored area, and thus increases the bulk power transfer through the monitored area. The contingencies are simulated under light, medium and heavy system loadings.

\begin{figure}[!htbp]
\vspace{-0.3cm}
 \centering
\subfloat[Contingency 1]{\includegraphics[width=45mm]{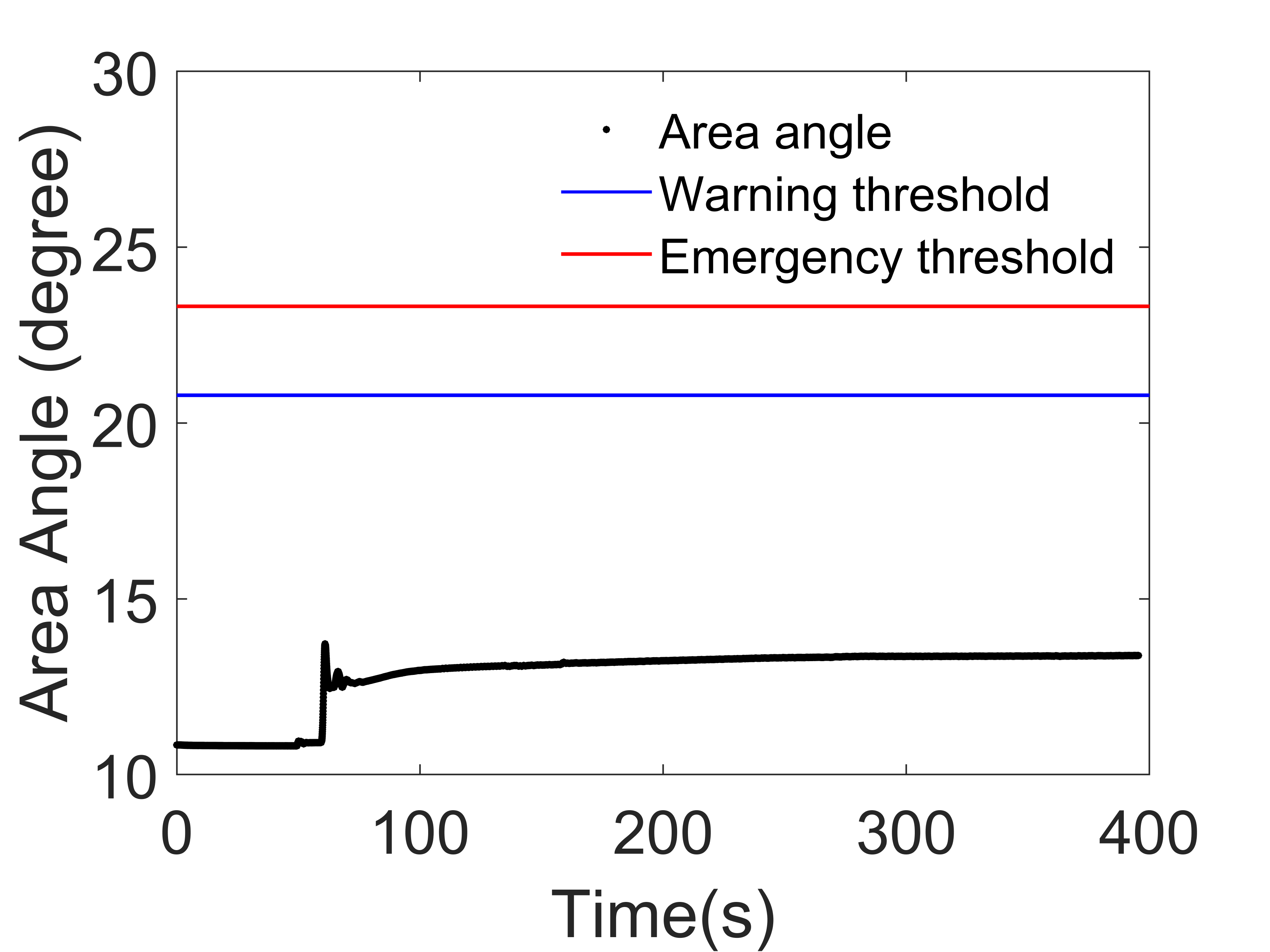}}
\subfloat[Contingency 2]{\includegraphics[width=45mm]{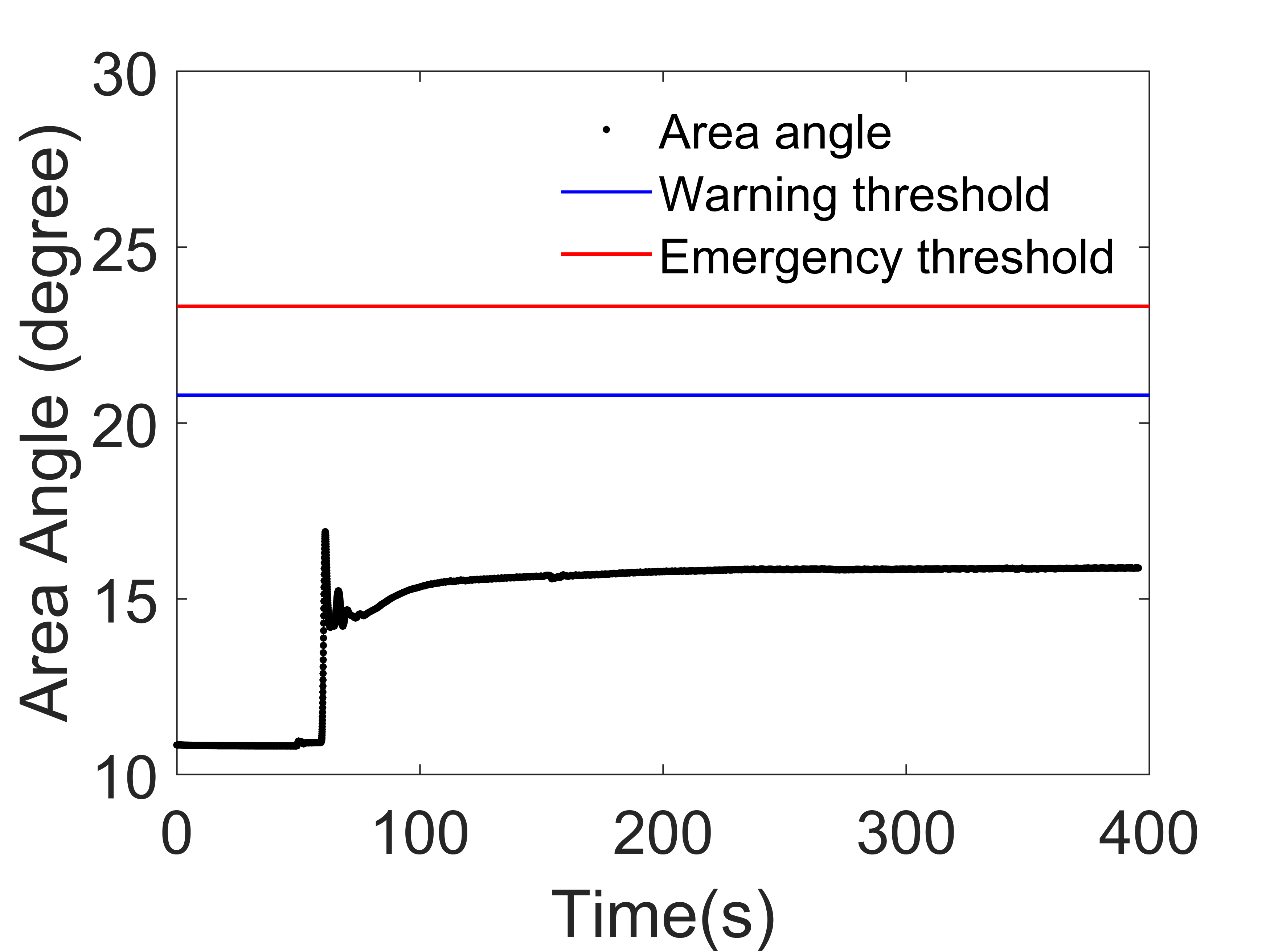}} 
\caption{Area angle under generator trip for light loading}
\vspace{-0.3cm}
\label{Fig9}
\end{figure}

\begin{figure}[!htbp]
\vspace{-0.5cm}
 \centering
\subfloat[Contingency 1]{\includegraphics[width=45mm]{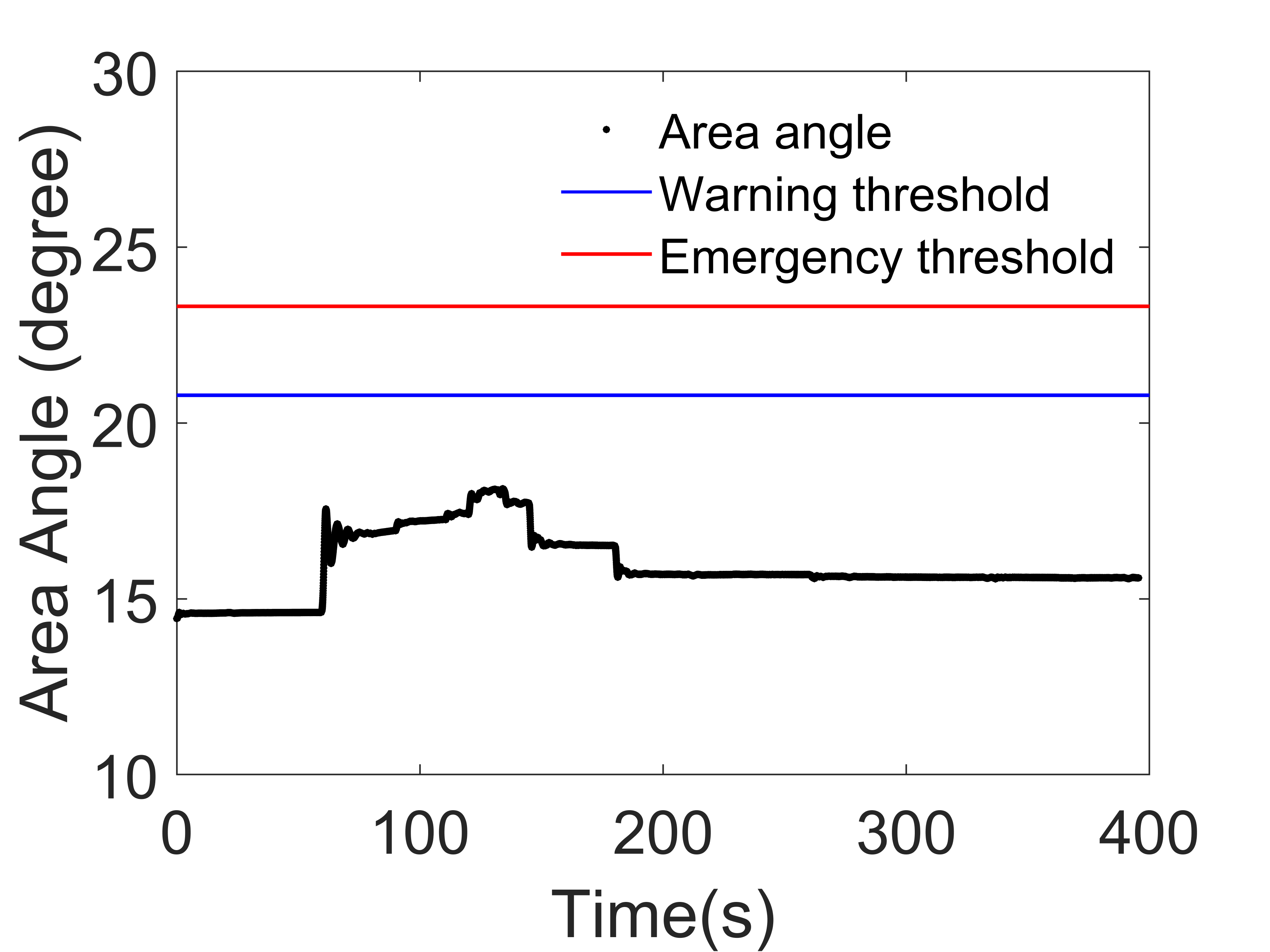}}
\subfloat[Contingency 2]{\includegraphics[width=45mm]{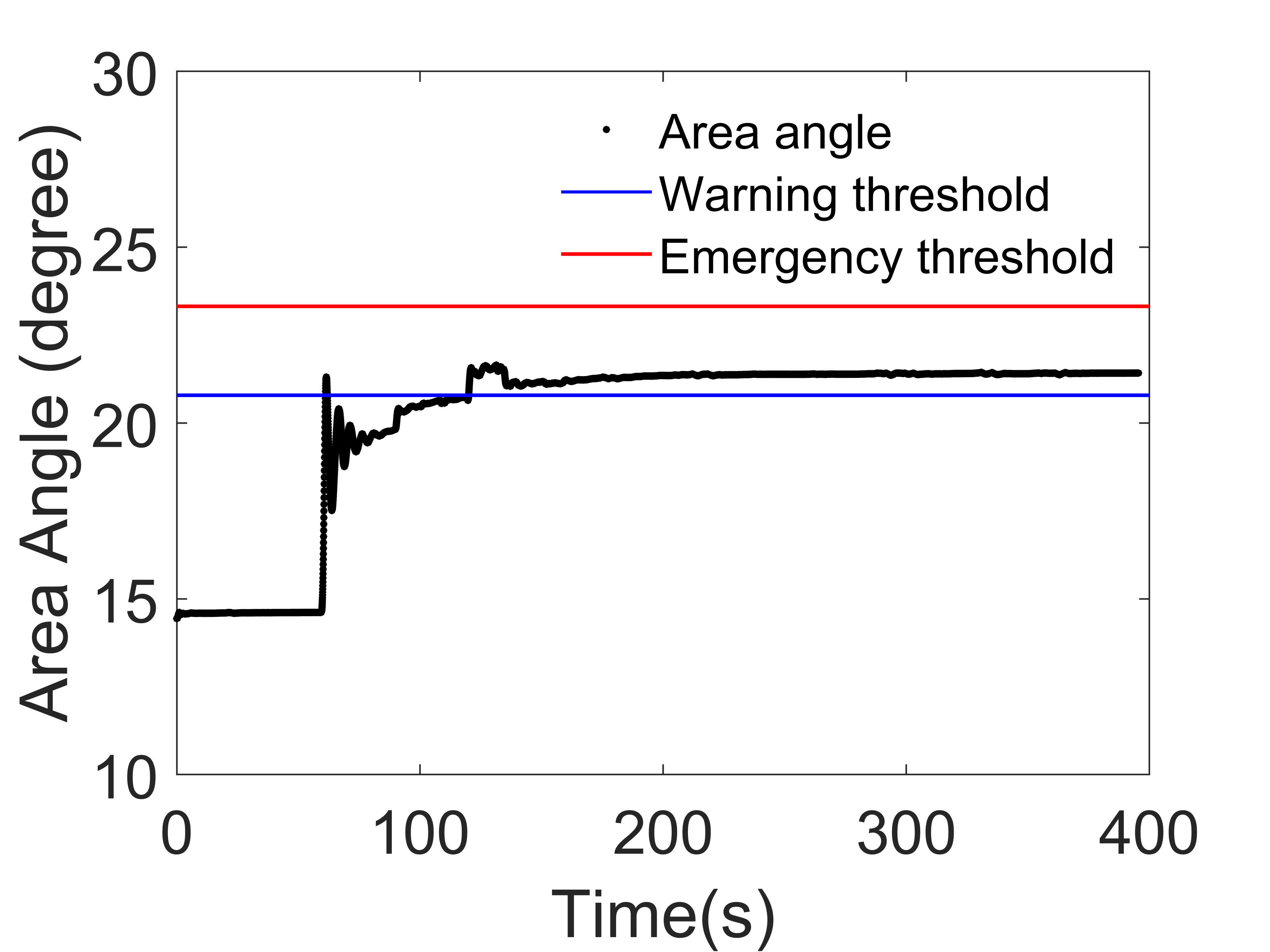}} 
\caption{Area angle under generator trip for medium loading}
\vspace{-0.3cm}
\label{Fig10}
\end{figure}

\begin{figure}[!htbp]
\vspace{-0.5cm}
 \centering
\subfloat[Contingency 1]{\includegraphics[width=45mm]{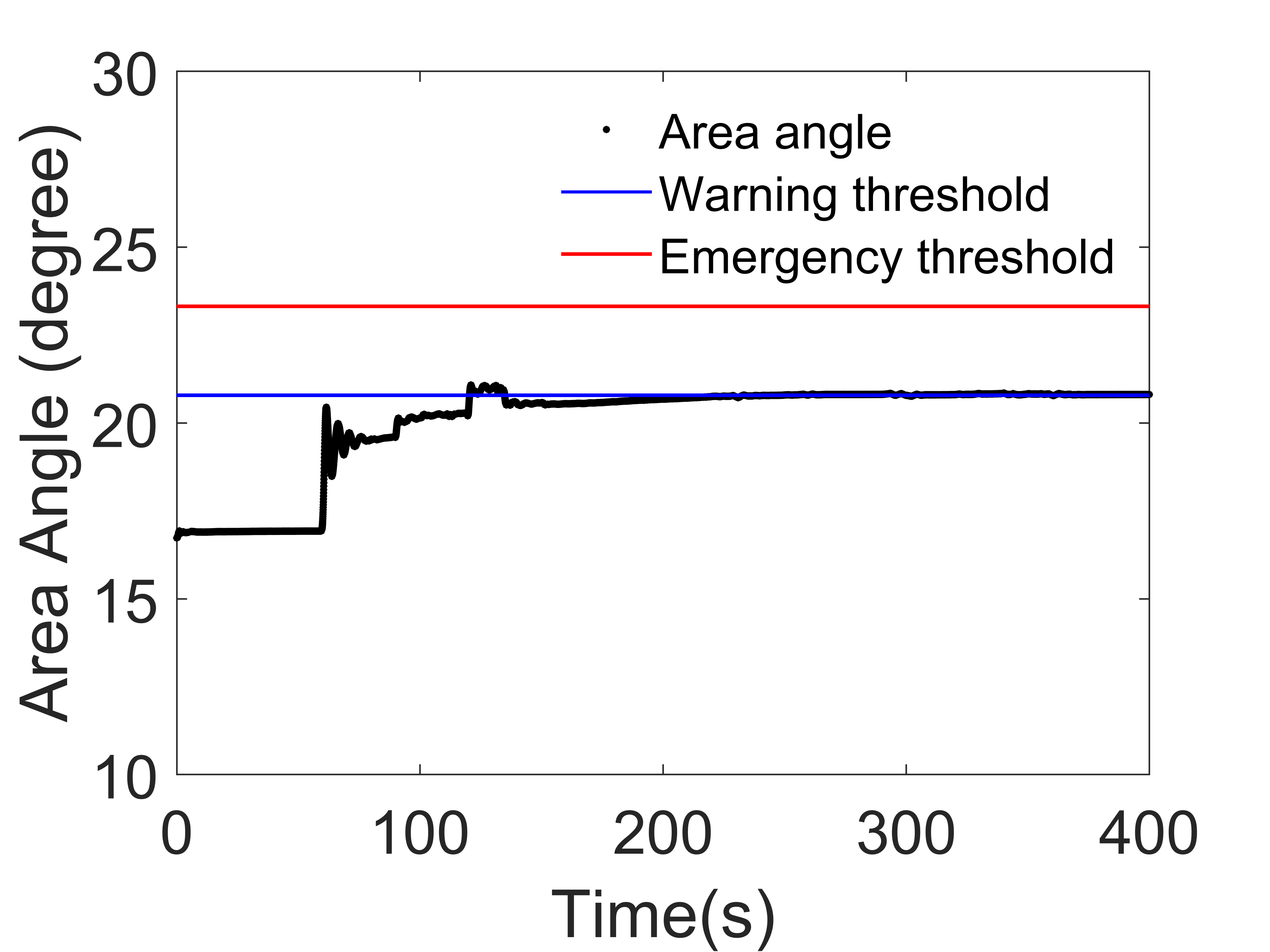}}
\subfloat[Contingency 2]{\includegraphics[width=45mm]{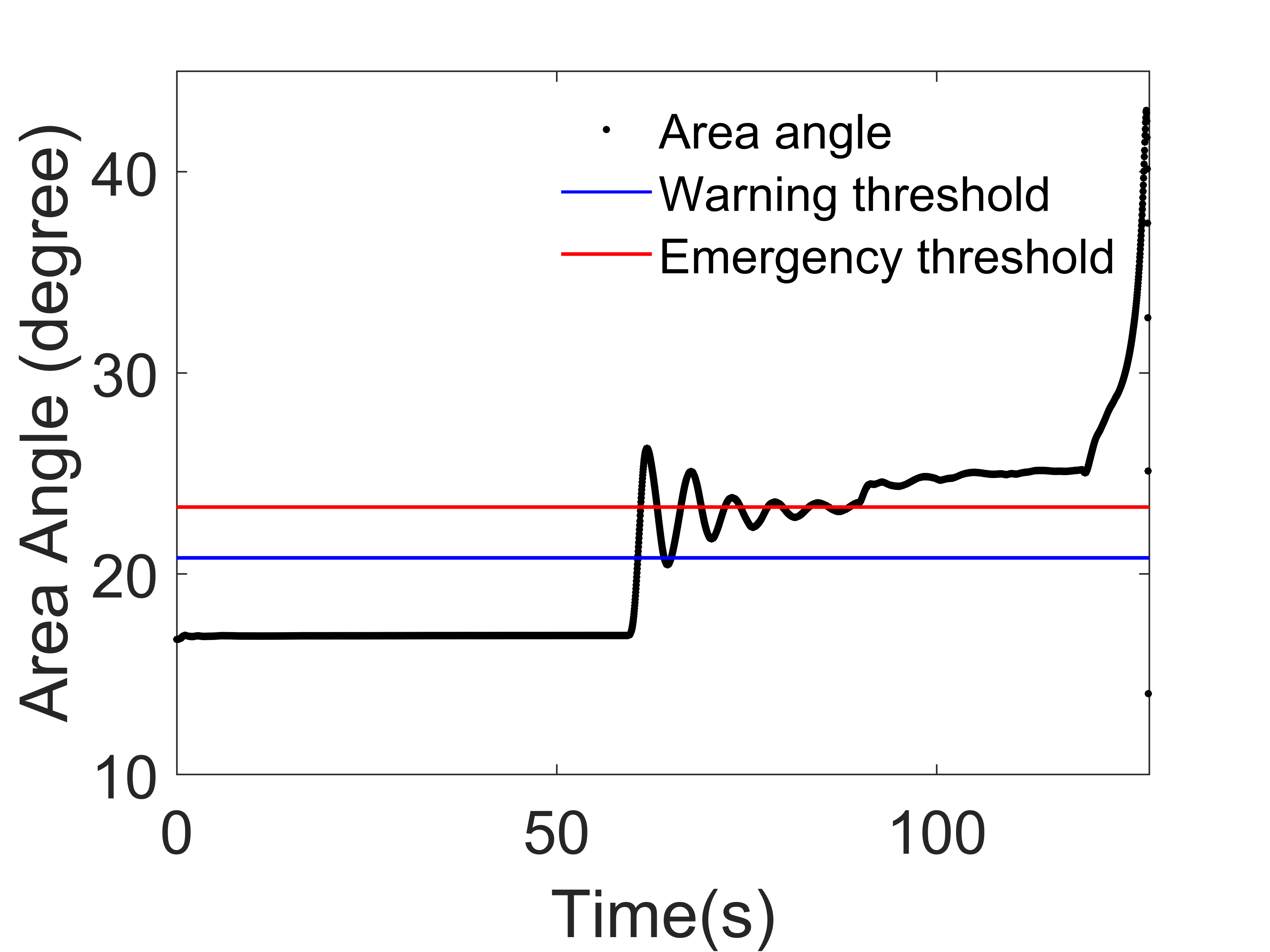}} 
\caption{Area angle under generator trip for heavy loading}
\label{Fig11}
\vspace{-0.3cm}
\end{figure}

From Fig.~\ref{Fig9}, we can see that the area angle increases after the contingency occurs, indicating the increased stress of bulk power transfer through the area. 
But the contingencies do not cause area angle to exceed any threshold.

From Fig.~\ref{Fig10}, we can see that the area angle increases after the contingency occurs. For Contingency 2, the area angle exceeds the warning threshold. The difference of area angle between Fig.~\ref{Fig9} and Fig.~\ref{Fig10} is caused by the increased power flow through the area as the system load level increases. 

From Fig.~\ref{Fig11}, we can see that Contingency 1 makes the area angle exceed the warning threshold and Contingency 2 makes the area angle exceed the emergency threshold. Warning status and emergency status will be indicated for Contingency~1 and Contingency 2, respectively. Note that the system goes unstable around 120\,s for Contingency 2. 

\textit{2) Verification of AAM with Line Outage}

\begin{figure}[!htbp]
\vspace{-0.6cm}
\centering
\subfloat[Contingency 3]{\includegraphics[width=45mm]{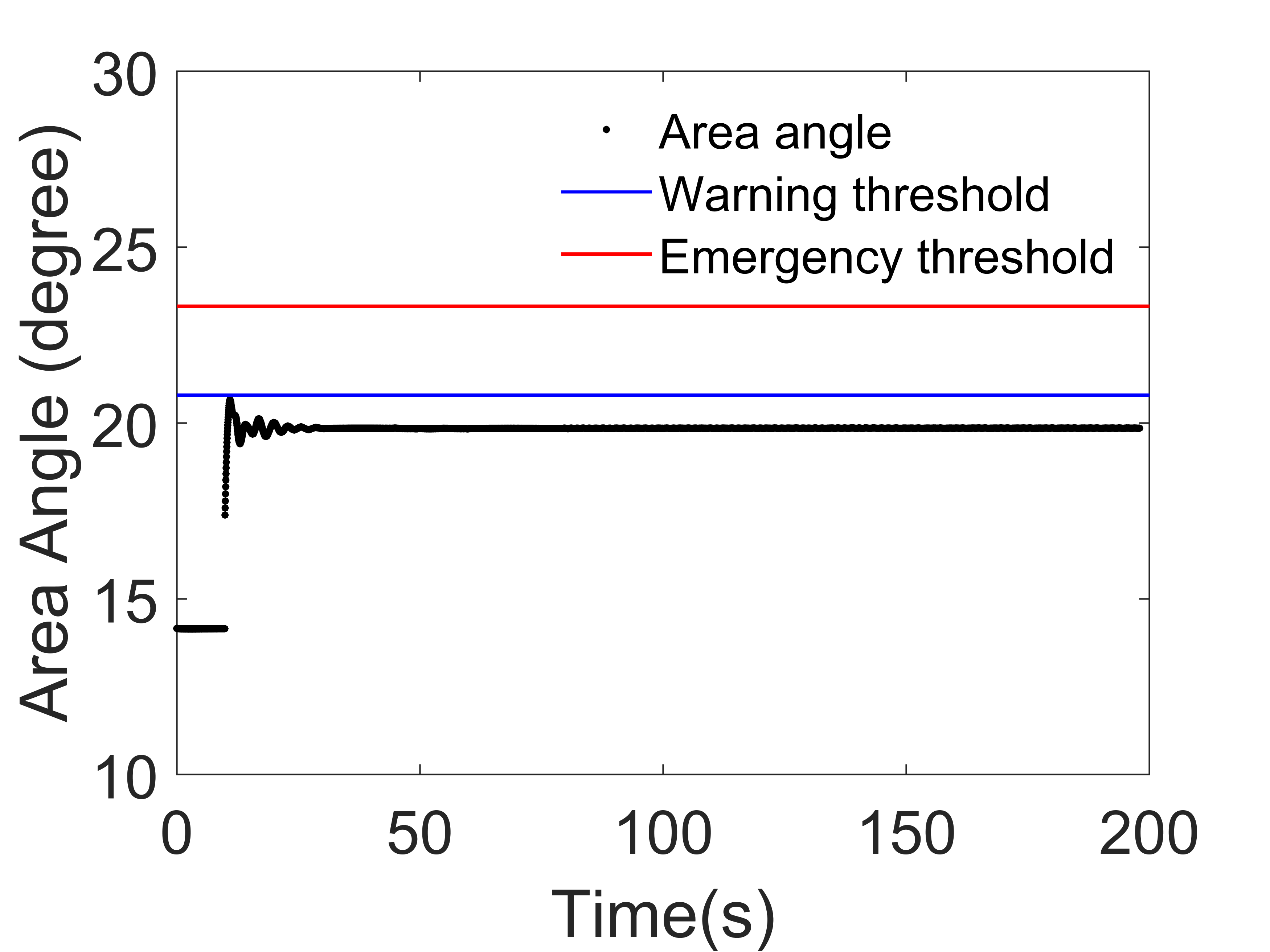}}
\subfloat[Contingency 4]{\includegraphics[width=45mm]{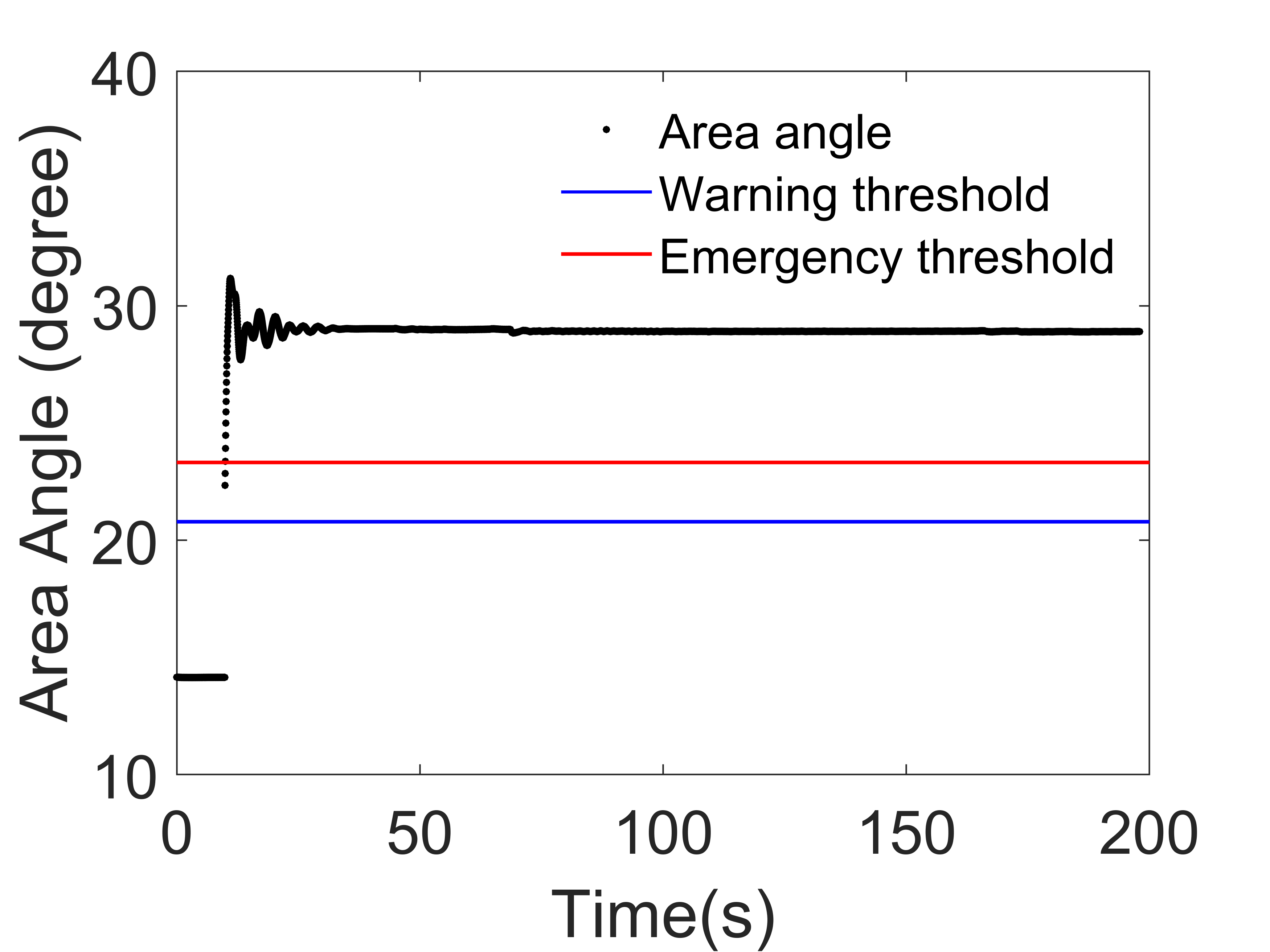}}
\caption{Area angle under line outage for medium loading}
\label{Fig12}
\vspace{-0.3cm}
\end{figure}

Contingency 3: Trip of one 500\,kV line at 10\,s.

Contingency 4: Trip of three 500\,kV lines at 10\,s.

These contingencies are inside the monitored area and they reduce the capability of bulk power transfer inside the monitored area and thus increase the stress. The contingencies are simulated under medium loading.

The area angle in Fig.~\ref{Fig12} (a) increases after this contingency occurs but it does not exceed any threshold.
Fig.~\ref{Fig12} (b) represents a more severe contingency and the area angle exceeds the emergency threshold. Emergency status will be indicated for Contingency 4.

\vspace{-0.1cm}
\subsection{Verification of AAM with Synchrophasor Measurements}

Two sets of recorded synchrophasor measurements from real contingencies are used here. These contingencies happened inside the monitored area. 

Contingency 5: Trip of one 500\,kV line.

Contingency 6: Trip of two 500\,kV lines. The time interval between the two line outages is around 100\,s. 

\begin{figure}[!htbp]
\vspace{-0.6cm}
\centering
\subfloat[Contingency 5]{\includegraphics[width=45mm]{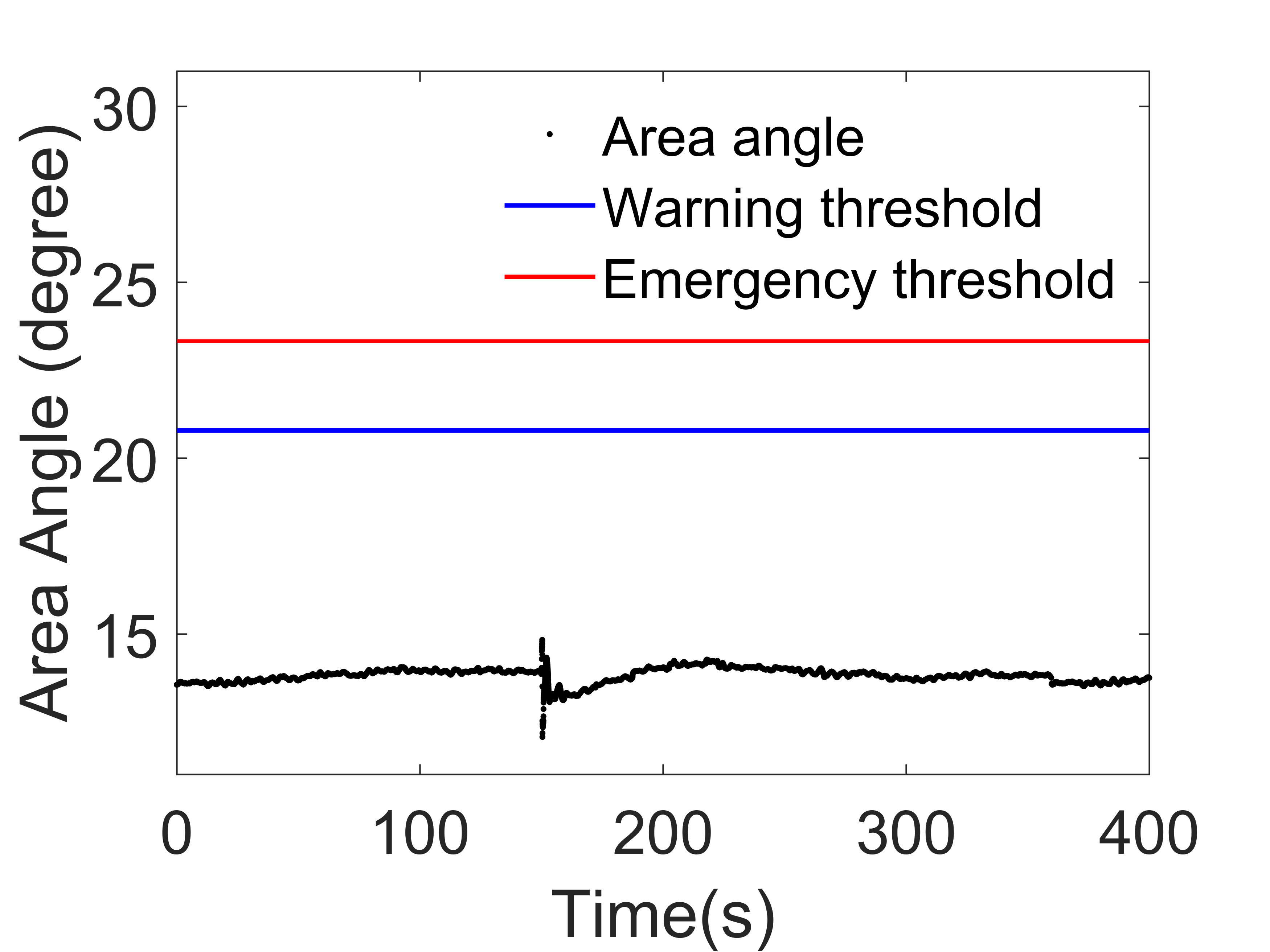}}
\subfloat[Contingency 6]{\includegraphics[width=45mm]{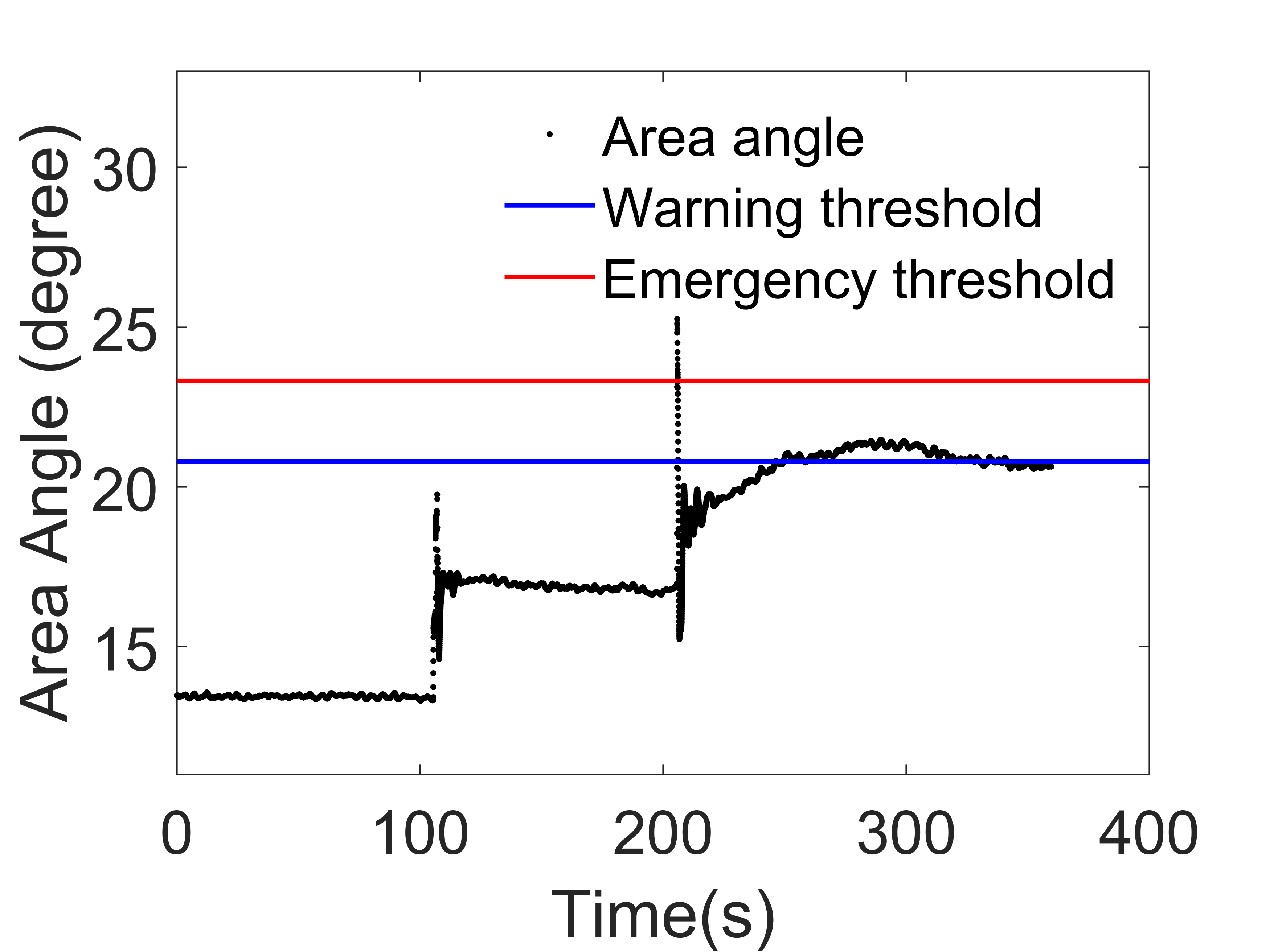}}
\caption{Area angle with synchrophasor measurements}
\label{Fig13}
\vspace{-0.3cm}
\end{figure}

From Fig.~\ref{Fig13} (a), we can see that the area angle varies for the whole time period but does not exceed any threshold.

From Fig.~\ref{Fig13} (b), we can see that the area angle increases significantly after the first line outage and continues increasing after the second line outage. It exceeds the warning threshold for more than 5\,s. The warning status will be indicated.

\vspace{-0.1cm}
\subsection{Influence of PACs on the Accuracy of Area Angle}
\label{PACinfluence}
The boundary buses [4, 5, 6, 7, 8, 9, 10, 14] do not have synchrophasor measurements. The PACs of them are calculated as [-1.8633, -3.1933, -1.0322, 2.7367, 0.0968, -8.7963, -9.5810, -3.9630] (degree). The weights are [0.0396, 0.0385, 0.005, 0.0023, -0.1269, -0.0958, -0.0017, -0.0395]. 



The influence of PACs on the accuracy of area angle is investigated with Contingency 1 under heavy loading.

\begin{figure}[!htbp]
\vspace{-0.5cm}
 \centering
\subfloat[Area angle comparison ]{\includegraphics[width=45mm]{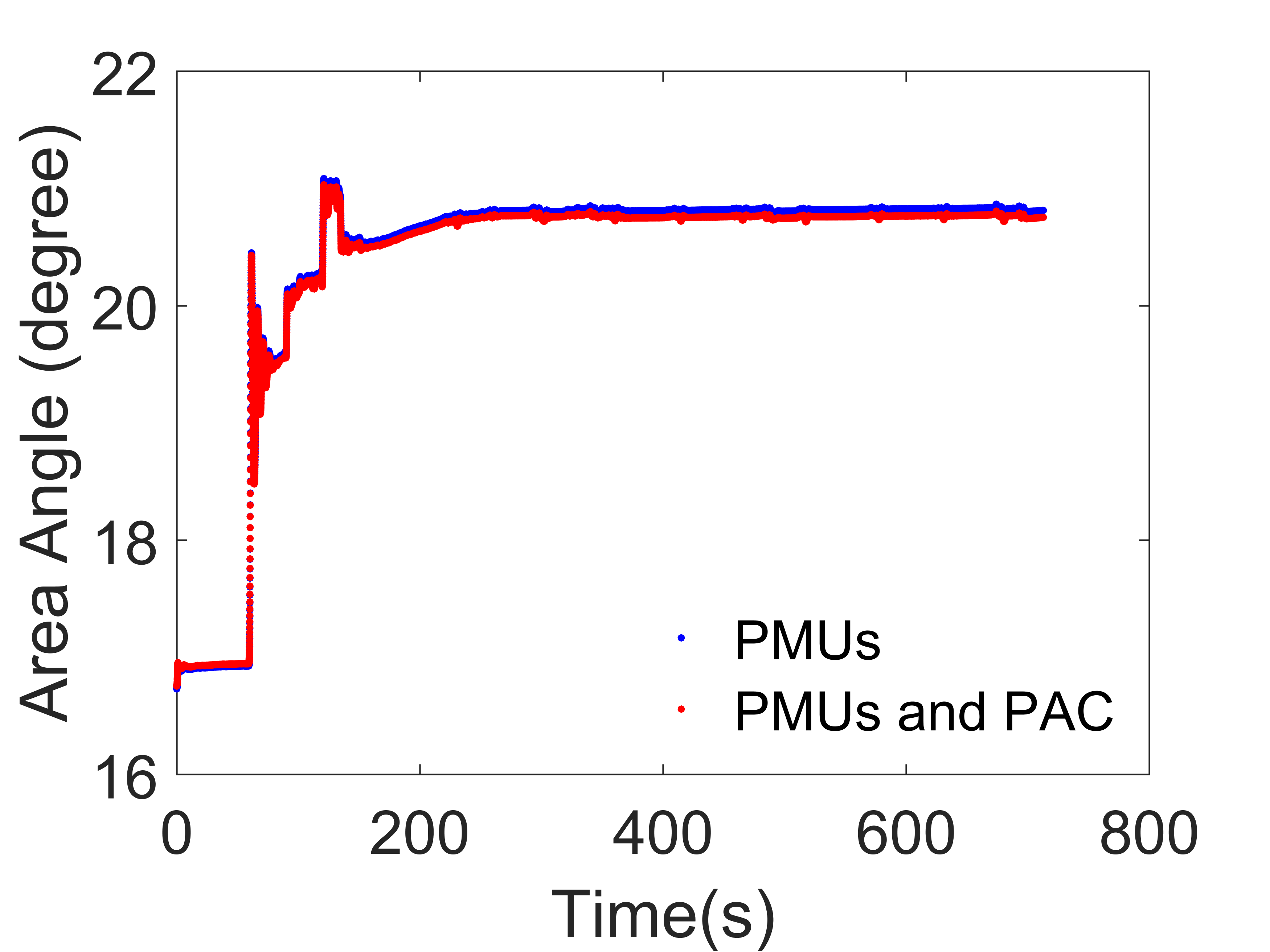}}
\subfloat[Area angle difference ]{\includegraphics[width=45mm]{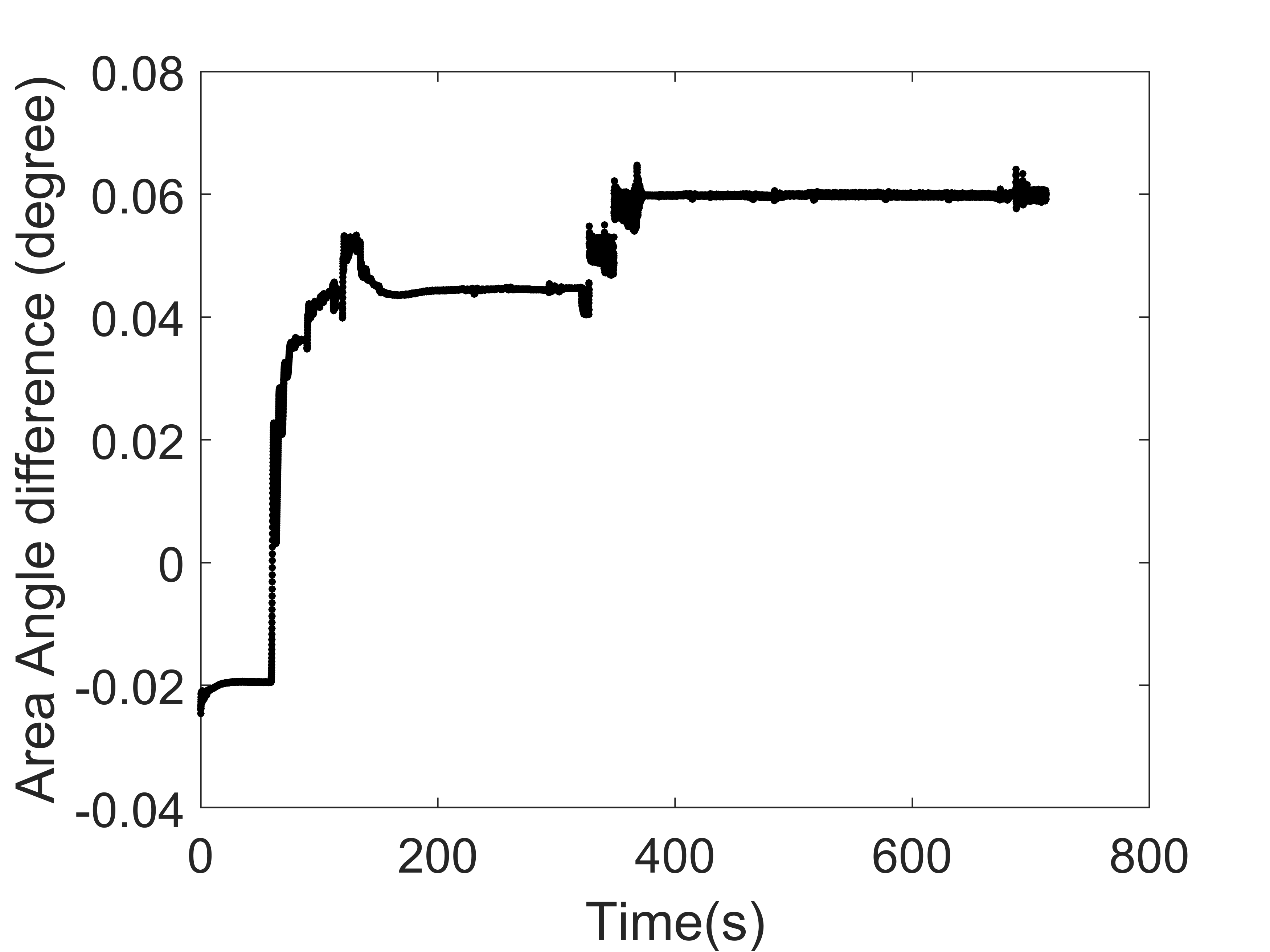}}
\caption{Influence of PACs on area angle for Contingency 1 }
\label{Fig14}
\vspace{-0.3cm}
\end{figure}

In Fig.~\ref{Fig14} (a), the curve marked by ``PMUs" assumes that all boundary buses are installed with synchrophasor measurements. The curve marked by ``PMUs and PAC" uses the Phase Angle Compensation method in Section III-E. From Fig.~\ref{Fig14} (b), we can see that the largest mismatch is around 0.06 degree. It suggests that the PAC method is accurate enough to calculate area angle. We can see that the weights for those boundary buses without synchrophasor measurements are quite small and the absolute values of PACs are not large, thus the influence on the accuracy of area angle is very small.

The PAC values obtained from a power flow model could vary with system operating condition. We investigate the effect of using different sets of PAC values on the area angle. Different sets of PAC values obtained from light, medium and heavy loading conditions with the reduced model are used with Contingency 1.
\begin{figure}[!htbp]
\vspace{-0.5cm}
\centering
\subfloat[PAC comparison]{\includegraphics[width=45mm]{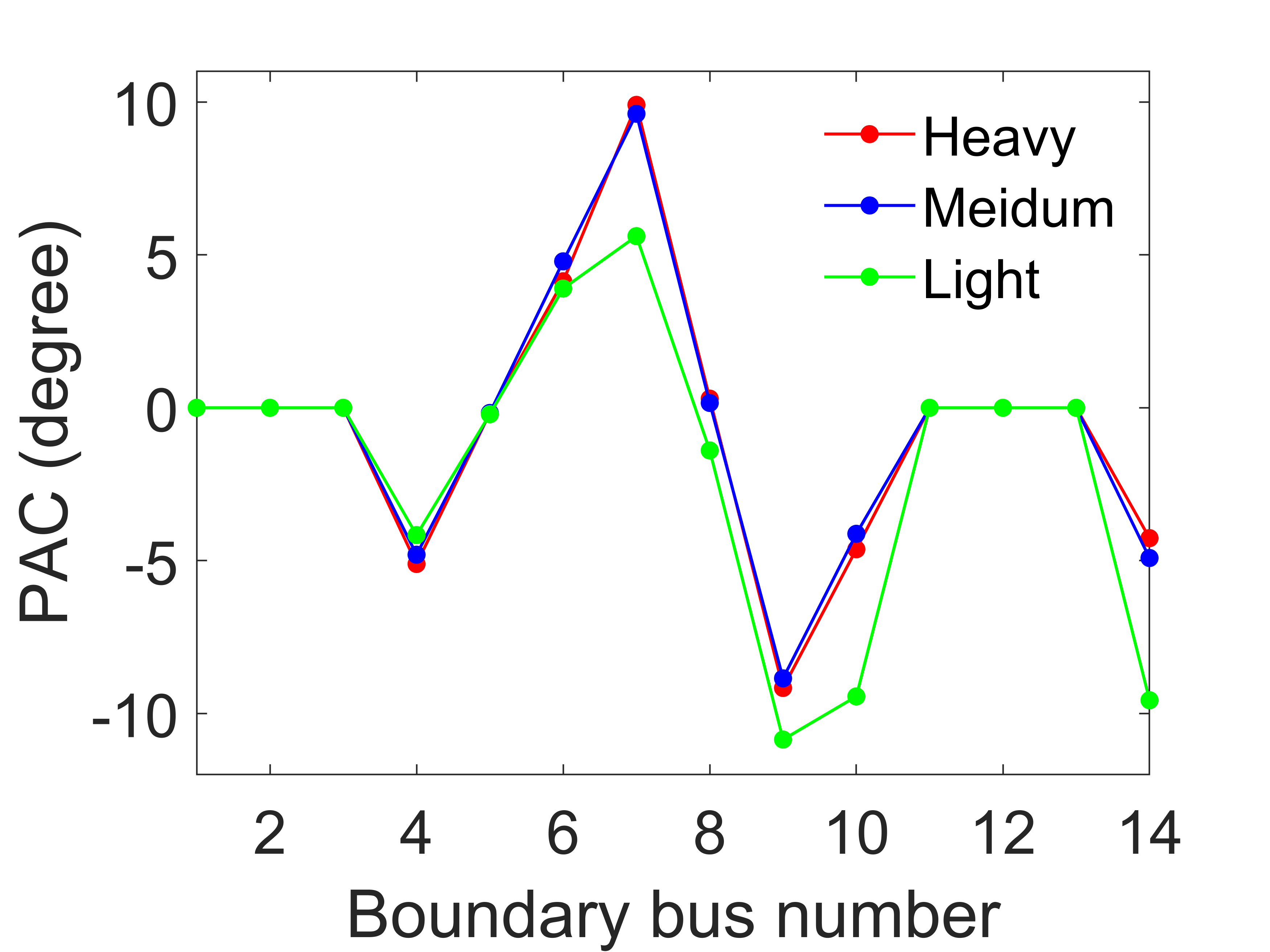}}
\subfloat[Area angle comparison]{\includegraphics[width=45mm]{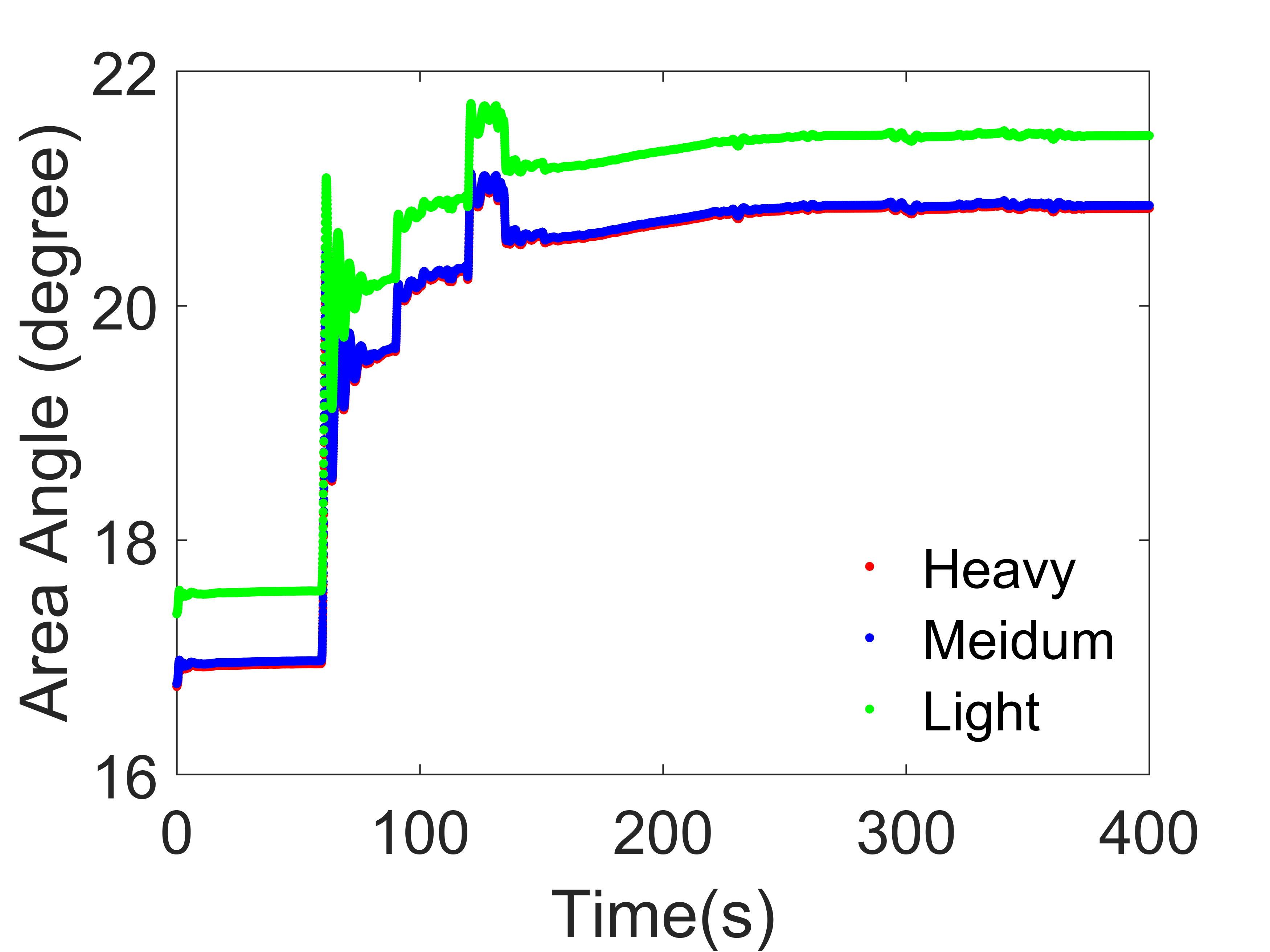}}
\caption{Influence of different sets of PACs on area angle}
\label{Fig16}
\vspace{-0.3cm}
\end{figure}
From Fig.~\ref{Fig16}, the area angle curves are almost overlapped for the sets of PACs obtained from the medium and heavy loadings. However, the area angle difference is relatively large for using two sets of PACs obtained from the medium and light loadings. This is caused by the large difference of loading level between medium and light loading conditions, which is approximately 12\% of the original system loading level. The larger the loading difference, the larger the difference of the PACs and thus the area angle. However, in real-time operations,  such a significant change of loading is rare. In summary, the accuracy for calculating area angle with a set of PACs is quite high for usual loading levels.

\vspace{-0.1cm}
\subsection{Updating Thresholds Under Significant Topology Change}
\label{testtopchange}

Consider the maintainance of two lines inside the monitored area, the updated warning and emergency thresholds using two methods are given in Table~\ref{tab4}. The results using the original method are benchmarks.

\begin{table}[!htbp]
\centering  
\vspace{-0.1cm}
\caption{Area Angle Thresholds}
\begin{tabular}{ccc}
 & Warning Threshold & Emergency Threshold  \\ \hline 
Original Method& $\theta^{thr,w}_{mod}$=22.08& $\theta^{thr,e}_{mod}$=26.75 \\ 
Proposed Method& $\hat{\theta}^{thr,w}_{mod}$=23.44 & $\hat{\theta}^{thr,e}_{mod}$=26.88 \\  \hline
\end{tabular}
\label{tab4}
\end{table}

From Table~\ref{tab4}, the mismatch of $\theta^{thr,w}_{mod}$ between two methods is 6.16\%, which is not small. The mismatch of $\theta^{thr,e}_{mod}$ between two methods is 0.49\%, which is very small. Since mitigation strategies are needed if the area angle exceeds the emergency threshold, we concern more about the accuracy of emergency threshold. Thus the accuracy of the proposed method for updating area angle thresholds is acceptable.

\vspace{-0.1cm}
\subsection{Mitigation Strategy for Reducing the Bulk Power Stress}
Contingencies 7, 8, and 9 are used to verify the proposed mitigation strategy in Section III-D.
Each contingency has three stages. The first stage is the normal state before the contingency. The second stage is immediately after the contingency. The third stage is  after the mitigation of load shedding on the buses of receiving side.

\begin{figure}[!htbp]
\vspace{-0.4cm}
\centering
\includegraphics[width=45mm]{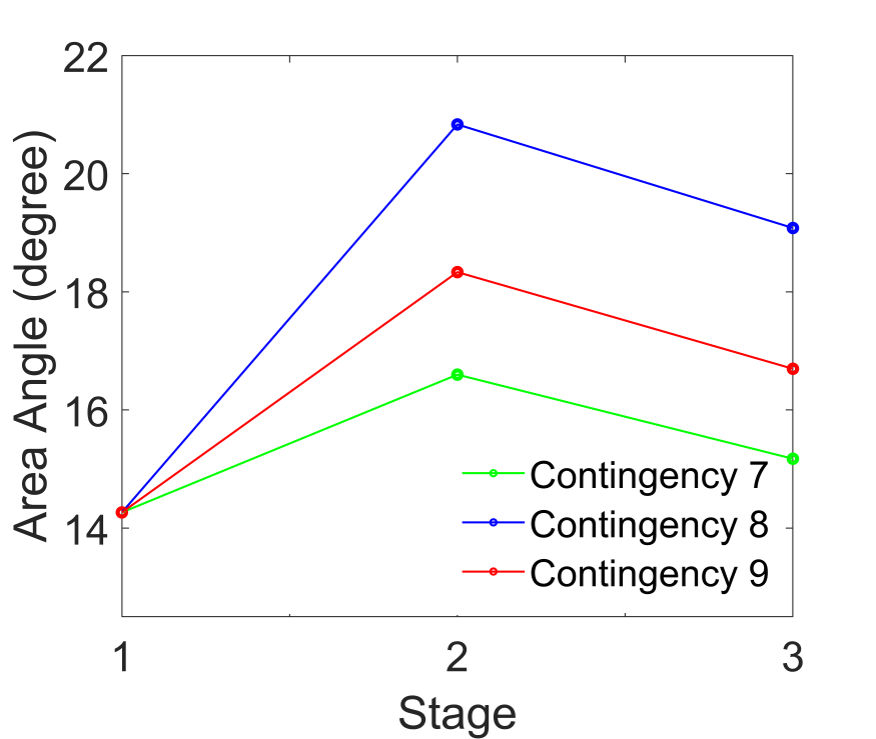}
\caption{Reduce area angle for three contingencies}
\label{Fig17}
\vspace{-0.2cm}
\end{figure}

From Fig.~\ref{Fig17}, we can see that the area angle increases at stage 2 compared with that of stage 1 and then decreases at stage 3 after performing load shedding compared with that of stage 2 for each contingency, which shows the proposed mitigation strategy working to reduce the bulk power stress.

\vspace{-0.1cm}
\subsection{Real-Time Application of AAM}
\begin{figure}[!htbp]
\centering
\includegraphics[width=82mm]{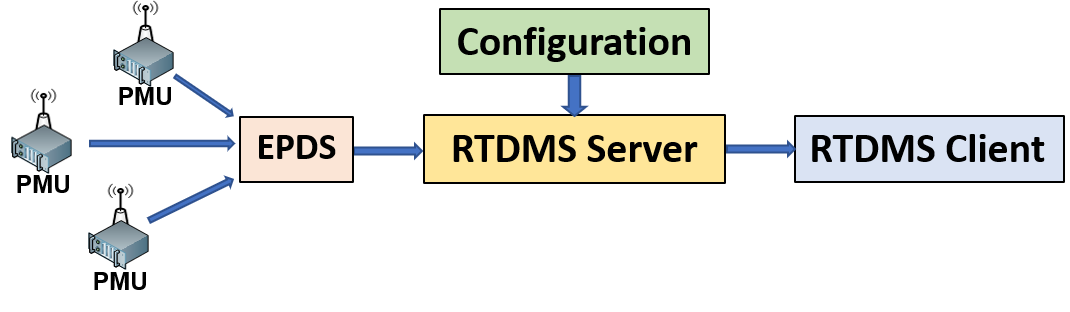}
\caption{Platform of real-time application of AAM}
\label{Fig18}
\vspace{-0.3cm}
\end{figure}

We
develop a platform for real-time application of AAM as shown by Fig.~\ref{Fig18}. Synchrophasor measurements collected by Enhanced Phase Data Simulator (EPDS) are sent out to Real Time Dynamics Monitoring System (RTDMS) Server through C37.118 data stream protocol. 
The area angle is calculated in the RTDMS Server. The area angle and area angle thresholds are visualized in the RTDMS Client in real-time.


This platform is deployed in BPA. It is running in their laboratory with live stream data from synchrophasor measurements. BPA has also tested  it using simulated data and recorded synchrophasor measurements for historical events.
\begin{figure}[!htbp]
\vspace{-0.3cm}
\centering
\includegraphics[width=84mm]{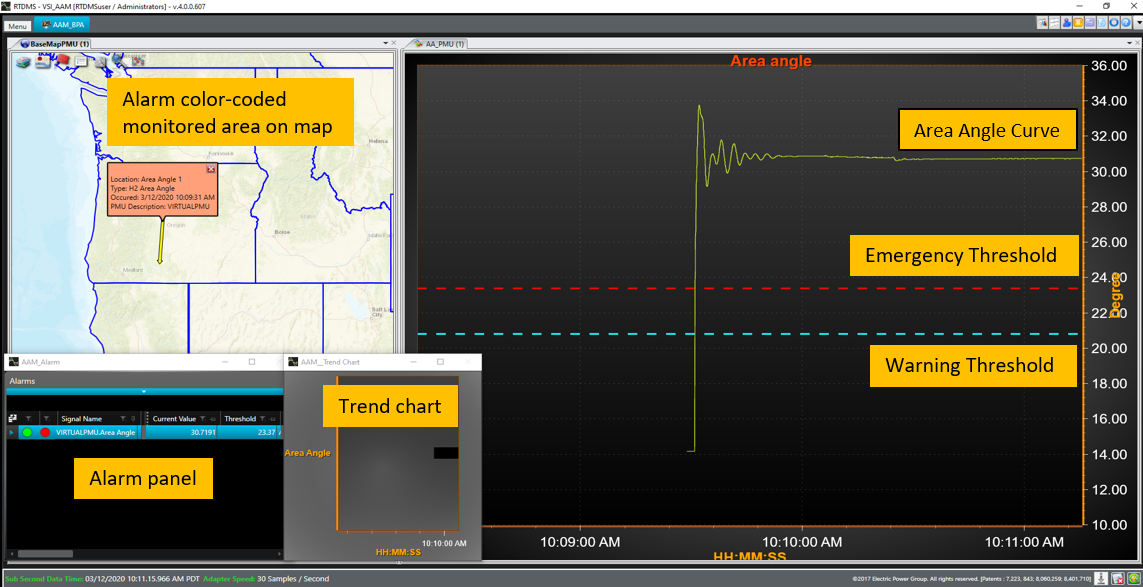}
\caption{Visualization of AAM in real-time}
\label{Fig19}
\vspace{-0.3cm}
\end{figure}
For example, the area angle can be seen in real-time responding to Contingency 4 under medium loading in Fig.~\ref{Fig19}. The emergency status is reported with red in the ``Alarm Panel" after the area angle exceeds the threshold. The mitigation strategy will be implemented into the platform in future work.

\section{Conclusion}
This paper develops and applies a practical framework of area angle monitoring (AAM) to monitor the stress of bulk power transfer across an area in real-time using synchrophasor measurements. To make AAM applicable in practice, we give methods to handle incomplete synchrophasor measurements at the boundary of the area, and methods to identify and quickly update the area angle thresholds for emergency and warning actions based on the AAM. We demonstrate the use of these mitigation actions.
Case studies and a utility deployment for an area demonstrate AAM with both simulated data and recorded and live-stream synchrophasor measurements. 
The innovations and testing of AAM position it as a practical tool for monitoring area stress and suggesting mitigation actions to the operators when thresholds are exceeded.


\begin{thebibliography}{11}
\bibitem{FA1}
L. Zhang, H. Chen, Q. Wang, N. Nayak, Y. Gong, A. Bose, ``A novel on-line substation instrument transformer health monitoring system using synchrophasor data," \textit{IEEE Trans. Power Delivery}. doi: 10.1109/TPWRD.2019.2905426.




\bibitem{FA2}
R. W. Cummings, in \textit{Predicting Cascading Failures, Presentation at NSF/EPRI Workshop on Understanding and Preventing Cascading Failures in Power Syst}., Westminster, CO, USA, Oct. 2005.

\bibitem{FA3}
V. Venkatasubramanian, Y. X. Yue, G. Liu, et al, ``Wide-area monitoring and control algorithms for large power systems using synchrophasors,"  \textit{IEEE Power Syst. Conf. and Exposition}, Seattle, WA, USA, Mar. 2019.

\bibitem{FA4}
U.S.-Canada Power System Outage Task Force, ``Final report on the August 14, 2003 blackout in the United States and Canada: Causes and Recommendations," Apr. 2004.







\bibitem{NERC02}  NERC (North American Electric Reliability Council),
 1996 system disturbances, 2002.


 
  \bibitem{AZFERCNERC12}  Arizona-Southern California Outages on September 8, 2011: Causes and Recommendations, Federal Energy Regulatory Commission and the North American Electric Reliability Corporation, April 2012.





\bibitem{FA9}
D. Gan, X. Luo, D. V. Bourcier, R. J. Thomas, ``Min-max transfer capability of transmission interfaces," \textit{Int. J. Elect. Power Energy Syst}., vol. 25, no. 5, pp. 347-353, 2003.

\bibitem{FA10}
Y. V. Makarov, P. Du, S. Lu, T. B. Nguyen, X. Guo, \textit{Wide Area Power System Security Region PNNL-19063}, Richland, WA USA:Pacific Northwest National Lab., 2009.

\bibitem{FA11}
F. Capitanescu, T. V. Cutsem, ``Evaluating bounds on voltage and thermal security margins under power transfer uncertainty," \textit{Proc. PSCC Conf}., Jun. 2002.

\bibitem{FA12}
Simultaneous Transfer Capability: Direction for Software Development, Electric Power Research Institute Report EL-7351 Project 3140-1, 1991.

\bibitem{FA13}
I. Dobson, ``Voltages across an area of a network," \textit{IEEE Trans. Power Syst}., vol. 27, no. 2, pp. 993-1002, May 2012.

\bibitem{FA14}
I. Dobson, M. Parashar, ``A cutset area concept for phasor monitoring," \textit{Proc. IEEE PES General Meeting}, Jul. 2010.

\bibitem{FA15}
A. Darvishi, I. Dobson, ``Threshold-based monitoring of multiple outages with PMU measurements of area angle," \textit{IEEE Trans. Power Syst}., vol. 31, no. 3, pp. 2116-2124, May. 2016.

\bibitem{FA60}
M. W, L. Xie, ``Online detection of low-quality synchrophasor measurements: a data-driven approach," \textit{IEEE Trans. Power Syst}, vol. 32, no. 4, pp. 2817-2827, Jul. 2017.

\bibitem{FA17}
L. Zhang, A. Bose, A. Jampala, et al, ``Design, testing, and implementation of a linear state estimator in a real power system," \textit{IEEE Trans. Smart Grid}, vol. 8, no. 4, pp. 1782-1789, Jul. 2017.

\bibitem{FA16}
D. Shi, D. L. Shawhan, N. Li, et al, ``Optimal generation investment planning: Pt. 1: network equivalents," \textit{North American Power Symposium}, Aug. 2012.

\bibitem{FA27}
F. Dorfler, F. Bullo, ``Kron reduction of graphs with applications to electrical network," \textit{IEEE Trans. Circuits and Systems-I}, vol. 60, no. 1, Jan. 2013.

\bibitem{FA18}
A. Darvishi, ``Monitoring of single and multiple line outages with synchrophasors in areas of the power system," PhD. thesis, Iowa State University, Ames, IA, 2015.

\bibitem{FA19}
A. Darvishi, I. Dobson, ``Synchrophasor monitoring of single line outages via area angle and susceptance," \textit{North American Power Symposium}, Pullman WA USA, Sept. 2014.



\end{thebibliography}
\end{document}